\documentclass[12pt,twoside]{amsart}
\usepackage{latexsym,amsfonts,amsmath,amssymb,soul}
\usepackage{enumerate,soul}
\usepackage[titletoc]{appendix}
\usepackage[usenames, dvipsnames]{color}
\numberwithin{equation}{section}

\usepackage[latin1]{inputenc}
\usepackage{amsthm}
\usepackage{enumerate}
\usepackage{amstext}
\usepackage{amssymb}
\usepackage{blindtext}

\newtheorem{defn}{Definition}[section]
\newtheorem{thm}[defn]{Theorem}
\newtheorem{corollario}[defn]{Corollary}
\newtheorem{lemma}[defn]{Lemma}
\newtheorem{prop}[defn]{Proposition}
\newtheorem{oss}[defn]{Remark}
\newenvironment{pf}[1][\unskip]{\medskip\noindent{\bf Proof #1.}\enspace}
{\hfill\newline\smallskip}

\theoremstyle{definition}

\catcode`\@=11

\renewcommand{\section}%
   {\setcounter{equation}{0}\@startsection {section}{1}{\z@}{-3.5ex plus -1ex
  minus -.2ex}{2.3ex plus .2ex}{\large\bf}}

\title[Regularity of global solutions...]{Regularity of global solutions of partial differential equations in non isotropic ultradifferentiable spaces via time-frequency methods}
\author{Claudio Mele, Alessandro Oliaro}
\address{Dipartimento di Matematica e Fisica ``E. De Giorgi'', Universit\`a del Salento,
C.P.193, I-73100 Lecce, Italy. E-mail address: claudio.mele1@unisalento.it}
\address{Dipartimento di Matematica ``G. Peano'', Universit\`a di Torino, Via Carlo Alberto, 10, I-10123 Torino, Italy. E-mail address: alessandro.oliaro@unito.it}

\begin{document}

\begin{abstract}
In this paper we study regularity of partial differential equations with polynomial coefficients in non isotropic Beurling spaces of ultradifferentiable functions of global type. We study the action of transformations of Gabor and Wigner type in such spaces and we prove that a suitable representation of Wigner type allows to prove regularity for classes of operators that do not have classical hypoellipticity properties. \\[0.1cm]
\end{abstract}

	\maketitle
	\thispagestyle{empty}

\section{Introduction}

In this paper we are concerned with the regularity of linear partial differential operators with polynomial coefficients in ultradifferentiable classes. The problem of regularity was first introduced by Shubin \cite{13} in the frame of Schwartz functions and tempered distributions; a linear operator $A :\mathcal{S}'\to\mathcal{S}'$ is said to be regular if the conditions $u\in\mathcal{S}'$, $Au \in\mathcal{S}$ imply that $u \in\mathcal{S}$. In \cite{13} Shubin formulates an hypoellipticity condition (in the global pseudodifferential calculus developed there), proving that such condition is sufficient to have regularity of the correponding operator. On the other hand, such hypoellipticity is far to be necessary, as there are several examples of operators which are not hypoelliptic but are regular; for instance, in \cite{Wong} the regularity of the Twisted Laplacian
\begin{equation}\label{TL}
L=\left( D_x+\frac{1}{2}M_y\right)^2+\left( D_y-\frac{1}{2}M_x\right)^2
\end{equation}
is proved, despite the fact that $L$ is not hypoelliptic in the sense of Shubin; in \cite{12}, a class of twisted operators containing the Twisted Laplacian is studied, and a characterization of regularity for twisted differential operators of second order is provided; also in this case, the twisted operators consider in \cite{12} are never hypoelliptic in the sense of Shubin. On the other hand, the problem of characterizing regularity for classes of operators is quite hard. Even in very particular cases (as for ordinary differential operators with polynomial coeffcients) necessary and sufficient conditions for regularity are not known; an interesting work in this sense is \cite{Nicola-Rodino:2}, where necessary and sufficient conditions for ordinary differential operators are found, but under additional conditions on the roots of their Weyl symbol. The notion of regularity (in global sense) can be defined each time we have a (global) space of functions and a corresponding space of (ultra)distributions, and operators acting on the ultradistribution space. In particular, this problem can be considered in the frame of ultradifferentiable classes. Ultradifferentiable spaces have been widely studied, starting from the work of Gevrey \cite{G}, who introduced an intermediate scales of spaces between real analytic and $C^\infty$ functions in order to analyze the regularity of solutions of partial differential equations. Then Komatsu \cite{K}, and later Beurling \cite{Beu} and Bj\"orck \cite{10} introduced a class of ultradifferentiable functions, where the regularity (and eventually the growth, when treating global spaces) is controlled by suitable weight sequences or suitable weight functions. In \cite{5}, general results on the spaces defined through weight functions are proved, and in \cite{BMM} a comparison between spaces defined by weight sequences and by weight functions is provided, proving that the two approaches have an intersections but there are spaces that can be defined only through one of them. Then, a large amount of papers on spaces of this kind have been produced, on properties of the spaces themselves, in connection with the behavior of solutions of partial differential equations, or more recently also in connection with time-frequency analysis. \\[0.1cm]
In this paper we study non isotropic spaces of global ultradifferentiable functions of Beurling type, following the approach of \cite{10,5}, and we analyze how tools from time-frequency analysis can be profitably used to find large classes of examples of partial differential equations that are regular in this ultradifferentiable setting. The connection between the partial differential equations world and the time-frequency analysis world has been profitably investigated in the last years, and has produced interesting results; we refer for instance to \cite{RW}, where the H\"ormander global wave front set has been re-defined through Gabor transform and Gabor frames, and to \cite{BJO-Gabor} and \cite{BJOS}, where tools from time-frequency analysis are used in order to study wave front set and nuclearity properties in the frame of ultradifferentiable spaces. In this paper we fix a collection of weight function $\omega_{j},\sigma_j$ for $j=1,\dots,N$ (see Definition \ref{weight}), and we denote by $\Omega$ and $\Sigma$ the functions
\begin{equation*}
\Omega=\omega_1\oplus \dots \oplus \omega_N, \quad \Sigma=\sigma_1\oplus\dots\oplus\sigma_N,
\end{equation*}
i.e., $\Omega(x):=\omega_1(x_1)+\dots+\omega_N(x_N)$, and analogously for $\Sigma$; then we consider the Fr\'echet space $\mathcal{S}^{\Sigma}_{\Omega}(\mathbb{R}^N)$, defined as the set of all functions $f\in L^1(\mathbb{R}^N)$ such that $f,\hat{f}\in C^\infty$ and
\begin{align*}
&\| \exp(\lambda\Omega)D^\alpha f\|_\infty <\infty,\; \text{for each} \; \lambda>0,\; \alpha\in\mathbb{N}^N_0,
 \\
&\|\exp(\lambda\Sigma)D^\alpha\hat{f}\|_\infty <\infty, \; \text{for each} \; \lambda>0,\; \alpha\in\mathbb{N}^N_0.
\end{align*}
Each of the weight can be non-quasianalytic or quasianalytic, as the results of the present paper hold in both cases; these spaces allow different behavior in different directions, as well as different decays of the function $f$ and of its Fourier transform $\hat{f}$. Observe that the Fourier transform is no longer an automorphism on $\mathcal{S}_\Omega^\Sigma$, as it is on the Bj\"orck space $\mathcal{S}_\omega$, but maps $\mathcal{S}^{\Sigma}_{\Omega}$ into $\mathcal{S}_{\Sigma}^{\Omega}$. Moreover, the spaces $\mathcal{S}_\Omega^\Sigma$ contain as particular cases Beurling spaces of Gelfand-Shilov type. We give different equivalent systems of seminorms for the space $\mathcal{S}^{\Sigma}_{\Omega}(\mathbb{R}^N)$, in the spirit of the results contained in \cite{8,9}, and we consider the problem of regularity in this setting. We say that an operator is $S_\Omega^\Sigma$-regular if the conditions $u\in(\mathcal{S}_\Omega^\Sigma)'$, $Au \in\mathcal{S}_\Omega^\Sigma$ imply that $u \in\mathcal{S}_\Omega^\Sigma$. We consider here partial differential operators with polynomial coefficients; our approach follows an idea that is already present in some works related to engineering applications, see \cite{Galleani-Cohen-1}, \cite{Galleani-Cohen-2}. In these papers some equations are analyzed, looking for the Wigner transform of the solution. Instead of finding  first a solution $u$, and then computing its Wigner transform, the equation itself is Wigner-transformed, and, in some cases, the new equation allow to find directly the the exact expression of the Wigner transform of the solution. This approach works well in the cases of partial differential equations with polynomial coefficients, and has been already used to study regularity properties of solutions of partial differential equations, see for instance \cite{8,12}, where regularity in classical Schwartz spaces and in isotropic ultradifferentiable classes is analyzed in dimension $1$. Here we study the action of transformations from time-frequency analysis, namely the Gabor transform
$$
V_g f(x,\xi)=\int_{\mathbb{R}^N} \exp(-it\xi)f(t)\overline{g(t-x)}\,dt
$$
and the following transform of Wigner type
\begin{equation*}
Wig[u](x,\xi)=\int_{\mathbb{R}^N}\exp(-it\xi)u\left(x+\frac{t}{2},x-\frac{t}{2}\right)\, dt,
\end{equation*}
on the spaces $\mathcal{S}_\Omega^\Sigma$, and we prove general results in arbitrary dimension that allow us to find large classes of partial differential operators that are $\mathcal{S}_\Omega^\Sigma$-regular. More precisely, we consider here operators $P$ in
$\mathbb{R}^{2N}$ with polynomial coefficients of the form
\begin{equation*}
P(x,y,D_x,D_y)=\sum_{|\alpha+\beta+\gamma+\mu|\leq m} c_{\alpha\beta\gamma\mu}x^\alpha y^\beta D_x^\gamma D_y^\mu,
\end{equation*}
with $x,y\in\mathbb{R}^N$, $c_{\alpha\beta\gamma\mu}\in\mathbb{C}$, and $m\in\mathbb{N}$. We show that $P$ is $\mathcal{S}^{\Sigma}_{\Omega}$-regular if and only if $\widetilde{P}$ is $\mathcal{S}_{\Omega_1}^{\Sigma_1}$-regular for suitable weights $\Omega_1$ and $\Sigma_1$, where $\widetilde{P}$ satisfies
\begin{equation*}
 Wig[Pu]=\widetilde{P}Wig[u].
\end{equation*}
This allows us to construct classes of operators that are regular in our ultradifferentiable setting. For instance, we prove that given a polynomial $p(z,\zeta)=\sum_{|\alpha+\beta|\leq m} c_{\alpha\beta}z^\alpha\zeta^\beta$ in $\mathbb{R}^{2N}$, with $c_{\alpha\beta}\in \mathbb{C}$ and $z,\zeta\in\mathbb{R}^N$, with $p(z,\zeta)\neq 0$ for every $(z,\zeta)\in \mathbb{R}^{2N}$, then the following operators are $\mathcal{S}^{\Sigma}_{\Omega}$-regular:
	\begin{eqnarray*}
		P_1 &=& \sum_{|\alpha+\beta|\leq m} c_{\alpha\beta} \left(\frac{x+y}{2}\right)^\alpha \left(\frac{D_x-D_y}{2}\right)^\beta, \\
		P_2 &=& \sum_{|\alpha+\beta|\leq m} c_{\alpha\beta} (D_x+D_y)^\alpha (y-x)^\beta, \\
		P_3 &=& \sum_{|\alpha+\beta|\leq m} c_{\alpha\beta} \left(x-\frac{D_y}{2}\right)^\alpha \left(x+\frac{D_y}{2}\right)^\beta,\\
		P_4 &=& \sum_{|\alpha+\beta|\leq m} c_{\alpha\beta} \left(y+\frac{D_x}{2}\right)^\alpha \left(\frac{D_x}{2}-y\right)^\beta.
	\end{eqnarray*}
Moreover, the Twisted Laplacian is $\mathcal{S}_{\omega_1\oplus\omega_2}^{\sigma_1\oplus\sigma_2}$-regular, for every weight functions $\omega_1$, $\omega_2$, $\sigma_1$, $\sigma_2$. We then prove similar results considering, instead of the transformation $Wig$, a general representation in the Cohen class, defined as $Q[u]=\kappa\star Wig[u]$, for a kernel $\kappa\in \mathcal{S}'(\mathbb{R}^{2N})$. \\[0.2cm]
The paper is organized as follows. Section 2 is devoted to the study of the space $\mathcal{S}^{\Sigma}_{\Omega}$ and is properties. In Sections 3 we analyze the action of the Gabor and Wigner transform on $\mathcal{S}^{\Sigma}_{\Omega}$ and in Sections 4 and 5 we study the $\mathcal{S}_\Omega^\Sigma$-regularity through Wigner-like transform and through Cohen class representations, giving some examples.

\section{Weight functions and the space $\mathcal{S}^{\Sigma}_{\Omega}(\mathbb{R}^N)$}

In this section we introduce the non-isotropic space $\mathcal{S}^{\Sigma}_{\Omega}(\mathbb{R}^N)$ of ultradifferentiable functions of Beurling type. We start with the definition of weight function in the sense of \cite{5}.

\begin{defn}\label{weight}
A continuous increasing function $\omega:[0,\infty)\to[0,\infty)$ is called a weight function if it satisfies the following properties:
\begin{itemize}
    \item[($\alpha$)] there exists $K\geq1$ such that $\omega(2x)\leq K(1+\omega(x))$ for every $x\geq0$;
	\item[($\beta$)] $\omega(x)=o(x)$ as $x \to \infty$;
	\item[($\gamma$)] there exist $a \in \mathbb{R}$, $b>0$ such that $\omega(x) \geq a+b\log(1+x)$, for every $x\geq0$;
	\item[($\delta$)] $\varphi_\omega(x)= \omega \circ \exp(x)$ is a convex function.
		\end{itemize}
	\end{defn}
Given a weight function we can extend $\omega: \mathbb{R} \to [0,\infty)$ by defining
$\omega(x)=\omega(|x|)$  for all $x \in \mathbb{R}$ (of course, in the same way we could extend $\omega $ to $\mathbb{R}^N$ for every $N$). The condition $(\beta)$ is weaker than the condition of non-quasianalyticity $\int_{1}^{\infty} \frac{\omega(t)}{1+t^2}\, dt < \infty$. When the latter condition is satisfied, the spaces that we are going to define shall contain non trivial compactly supported functions. All the results of this paper hold under condition $(\beta)$, i.e., both in the non-quasianalytic and in the quasianalytic case.
\\[0.1cm]
As standard, we define the \emph{Young conjugate} $\varphi^*_\omega$ of $\varphi_\omega$ as
\begin{equation}\label{Yconj}
\varphi^*_\omega(s):=\sup_{t\geq 0}\{ ts-\varphi_\omega(t)\},\quad s\geq 0.
\end{equation}
We recall that $\varphi^*_\omega$ is an increasing convex function on $[0,+\infty)$ and it satisfies $\varphi_\omega^{**}=\varphi_\omega$.
\\[0.1cm]
We now recall some known facts that shall be useful in the following. At first, there is no loss of generality in assuming that $\omega|_{[0,1]}\equiv0$; as a consequence, we easily have from \eqref{Yconj} that $\varphi_\omega^*(0)=0$. Moreover the properties in the next proposition hold; they are well-known and can be found in many references, we refer for instance to \cite{9}, where (in Section 2 and in the Appendix) several basic properties of weights are collected and proved with minimal assumptions.
\begin{prop}\label{wheight prop gamma'}
	Let $\omega$ be a weight function. Then
	\begin{itemize}
		\item[(1)] For each $\lambda>0$, $j\in \mathbb{N}_0$ and $x\geq 0$ we have
		\begin{equation*}
		x^j \exp(-\lambda\omega(x))\leq \exp\left(\lambda\varphi^*\left(\frac{j}{\lambda}\right)\right);
		\end{equation*}
		\item[(2)] For each $\lambda>0$ and $x\geq1$ we have
		\begin{equation*}
		\underset{j\in \mathbb{N}}{\inf}\; x^{-j}\exp\left(\lambda\varphi^*\left(\frac{j}{\lambda}\right)\right)\leq \exp\left(-\left(\lambda-\frac{1}{b}\right)\omega(x)-\frac{a}{b}\right),
		\end{equation*}
		where $a,b$ are the constants appearing in condition $(\gamma)$;
		\item[(3)] There exists a constant $L>0$, depending on $\omega$, such that for every $\lambda>0$, $\rho \geq 1$ and $j\in \mathbb{N}_0$,
		\begin{equation*}
		\rho^j \exp\left(\lambda\varphi^*\left(\frac{j}{\lambda}\right)\right)\leq C_{\rho,\lambda}\exp\left(\lambda'\varphi^*\left(\frac{j}{\lambda'}\right)\right),
		\end{equation*}
        for each $0\leq \lambda'\leq \frac{\lambda}{L^{[\log \rho +1]}}$ and
        for a suitable constant $C_{\rho,\lambda}>0$;

		\item[(4)] For each $\lambda>0$ and $j\in\mathbb{N}_0$

		\begin{equation*}
		j!\leq C_\lambda \exp\left(\lambda\varphi^*\left(\frac{j}{\lambda}\right)\right),
		\end{equation*}
		for a suitable constant $C_\lambda>0$;

		\item[(5)] For each $x,y\geq 0$
		\begin{equation*}
		\omega(x+y)\leq K(1+\omega(x)+\omega(y)),
		\end{equation*}
		where $K$ is the constant appearing in condition $(\alpha)$. Observe that this condition is weaker than subadditivity (i.e. $\omega(x+y)\leq \omega(x) + \omega(y)$). The weight functions satisfying ($\alpha$) are not necessarily subadditive in general.
		\end{itemize}
\end{prop}
Consider a collection of weight function $\omega_{j},\sigma_j$ for $j=1,\dots,N$. We denote by $\Omega$ and $\Sigma$ the functions on $\mathbb{R}^N$ defined by
\begin{equation}\label{directSum}
	\Omega=\omega_1\oplus \dots \oplus \omega_N, \quad \Sigma=\sigma_1\oplus\dots\oplus\sigma_N,
\end{equation}
in the sense that $\Omega(x):=\omega_1(x_1)+\dots+\omega_N(x_N)$, and analogously for $\Sigma$, for $x\in\mathbb{R}^N$. We can suppose without loss of generality that all the $\omega_j$ and all the $\sigma_j$ satisfy condition $(\alpha)$ of Definition \ref{weight} with the same constant $K$. Similarly we assume that $\omega_j$ and $\sigma_j$ satisfy condition $(\gamma)$ with the same constants $a$ and $b$, for every $j=1,\dots,N$.
\\[0.1cm]
We define the following space of rapidly decreasing ultradifferentiable functions of Beurling type.
\begin{defn}
	Let $\Omega$, $\Sigma$ be weight functions as in \eqref{directSum}. We define $\mathcal{S}^{\Sigma}_{\Omega}(\mathbb{R}^N)$ as the space of all functions $f\in L^1(\mathbb{R}^N)$ such that $f,\hat{f}\in C^\infty(\mathbb{R}^N)$ and satisfy
	\begin{align}
	&\| \exp(\lambda\Omega)D^\alpha f\|_\infty <\infty,\; \text{for each} \; \lambda>0,\; \alpha\in\mathbb{N}^N_0,
\label{cond 1 swo} \\
	&\|\exp(\lambda\Sigma)D^\alpha\hat{f}\|_\infty <\infty, \; \text{for each} \; \lambda>0,\; \alpha\in\mathbb{N}^N_0.
\label{cond 2 swo}
	\end{align}
\end{defn}
The corresponding (countable) family of seminorms
$$
\| \exp(n\Omega) D^\alpha f\|_\infty,\quad \|\exp(m\Sigma) D^\beta \hat{f}\|_\infty,
$$
with $n,m\in\mathbb{N}$, $\alpha,\beta\in\mathbb{N}^N_0$, induces a topology of Fr\'echet space on $\mathcal{S}^{\Sigma}_{\Omega}(\mathbb{R}^N)$.
\begin{oss}\label{diffweight}
The weight functions $\omega_j$ and $\sigma_j$ in \eqref{directSum} of course do not need to be different. In the case when some of the $\omega_j$ (or some of the $\sigma_j$) coincide we can put together the corresponding variables, in the following sense: if for instance $\omega_2=\omega_1$, we can choose
$$
\Omega=\omega_1(x_1)+\omega_1(x_2)+\omega_3(x_3)+\dots+\omega_N(x_N),
$$
as in \eqref{directSum}, or also
$$
\Omega=\omega_1(|(x_1,x_2)|)+\omega_3(x_3)+\dots+\omega_N(x_N),
$$
and the corresponding space $\mathcal{S}_\Omega^\Sigma(\mathbb{R}^N)$ does not change. We have indeed that
\begin{equation*}
\begin{split}
\omega_1(x_1)+\omega_1(x_2) &\leq 2\omega_1(|(x_1,x_2)|)\leq 2\omega_1(|x_1|+|x_2|)\leq 2K(1+\omega_1(x_1)+\omega_1(x_2)),
\end{split}
\end{equation*}
since $\omega_1$ is increasing and satisfies Proposition \ref{wheight prop gamma'} (5). In particular, if $\omega_1=\dots =\omega_N=\sigma_1=\dots =\sigma_N:=\omega$, the space $\mathcal{S}_\Omega^\Sigma(\mathbb{R}^N)$ coincides with the space $\mathcal{S}_\omega$ considered for instance in \cite{10}, \cite{8}, \cite{9}, \cite{7}.
\end{oss}

\begin{oss}\label{in s}
	The condition $(\gamma)$ ensures us that, for $\Omega$ and $\Sigma$ as in \eqref{directSum}, the space $\mathcal{S}^{\Sigma}_\Omega(\mathbb{R}^N)$ is contained in $\mathcal{S}(\mathbb{R}^N)$, with continuous inclusion. Indeed, since $\log |x|\leq \sum_{j=1}^N \log(1+|x_j|)$, given $f \in \mathcal{S}^{\Sigma}_\Omega(\mathbb{R}^N) $ and $\alpha,\beta\in\mathbb{N}^N_0$, we have that
	\begin{align*}
		\|x^\alpha D^\beta f\|_\infty\leq \sup_{x \in \mathbb{R}^N} \exp(|\alpha|\log|x|)|D^\beta f(x)| \leq C\sup_{x \in \mathbb{R}^N}\exp\left(\frac{|\alpha|}{b}\Omega(x)\right)|D^\beta f(x)| <\infty,
	\end{align*}
	for some constant $C$, where $b$ is the constant appearing in condition $(\gamma)$, common for all $\omega_j$.
	\\[0.1cm]
Then we can rewrite the definition of $\mathcal{S}^{\Sigma}_\Omega(\mathbb{R}^N)$ as the set of all the rapidly decreasing functions that satisfy $(\ref{cond 1 swo})$ and $(\ref{cond 2 swo})$.
\end{oss}

\begin{oss}\label{s}
The following facts can be easily proved.
\begin{itemize}
\item[(a)]
	Given $f \in \mathcal{S}(\mathbb{R}^N)$, since $\mathcal{F}\mathcal{F}(f)=(2\pi)^{N}Rf$, where $Rf(x)=f(-x)$ is the reflection operator, we have that $f \in \mathcal{S}^{\Sigma}_{\Omega}(\mathbb{R}^N)$ if and only if $\hat{f} \in \mathcal{S}^{\Omega}_{\Sigma}(\mathbb{R}^N)$. Moreover the Fourier transform $\mathcal{F}:\mathcal{S}^{\Sigma}_{\Omega}(\mathbb{R}^N)\to \mathcal{S}^{\Omega}_{\Sigma}(\mathbb{R}^N)$ is a continuous isomorphism.
\item[(b)]
The space $\mathcal{S}_\Omega^\Sigma (\mathbb{R}^N)$ is closed under convolution, arithmetic product of functions, translation and modulation, where the translation and modulation operators are defined by $T_s f(x):=f(x-s)$ and $M_t f(x):=e^{itx} f(x)$, respectively, where $s,t,x\in\mathbb{R}^N$.
\end{itemize}
\end{oss}
\begin{defn}\label{def2weightS}
	Let $\Omega$, $\Sigma$ be weight functions as in \eqref{directSum}. We define $(\mathcal{S}_{\Omega}^{\Sigma})'(\mathbb{R}^N)$ as the set of the linear and continuous maps from $\mathcal{S}^{\Sigma}_{\Omega}(\mathbb{R}^N)$ to $\mathbb{C}$.
\end{defn}
The next two lemmas are proved in \cite{7}, in the case of subadditive weight functions; the proof in our case is strictly analogous and is omitted.
\begin{lemma}\label{prop swo}
	Let $\Omega$, $\Sigma$ be weight functions as in \eqref{directSum} and consider  $f,g \in \mathcal{S}^{\Sigma}_{\Omega}(\mathbb{R}^N)$, $\lambda>0$ sufficiently large. Then
	\begin{align*}
		&\| \exp(\lambda\Omega)(f\star g)\|_\infty \leq C_\lambda \| \exp(K\lambda\Omega)f\|_\infty\| \exp(K\lambda\Omega) g\|_\infty,\\&
		\| \exp(\lambda\Sigma)(\hat{f}\star \hat{g})\|_\infty \leq C_\lambda \| \exp(K\lambda\Sigma)\hat{f}\|_\infty\| \exp(K\lambda\Sigma) \hat{g}\|_\infty,
	\end{align*}
	for a suitable positive constant $C_\lambda$.
\end{lemma}
\begin{lemma}\label{lemma gz}
	Let $\Omega$, $\Sigma$ be weight functions as in \eqref{directSum} and consider  $f,g \in \mathcal{S}^\Sigma_\Omega(\mathbb{R}^N)$, $\lambda>0$. Then the following properties hold:
	\begin{itemize}
		\item[(1)] For every $\alpha\in\mathbb{N}^N_0$,
		\begin{equation*}
		\exp(\lambda\Omega(t))D^\alpha(M_\xi T_xg)(t)= \sum_{\beta \leq \alpha} \binom{\alpha}{\beta} (i\xi)^\beta M_\xi T_x \left(\exp(\lambda\Omega(x+\cdot))D^{\alpha-\beta}g(\cdot)\right)(t);
		\end{equation*}
		\item[(2)] For every $\alpha\in\mathbb{N}^N_0$
		\begin{equation*}
		\|\exp(\lambda\Omega)D^\alpha(M_\xi T_xg)\|_\infty \leq C\exp(K\lambda\Omega(x)) \sum_{\beta \leq \alpha} \binom{\alpha}{\beta} |\xi^\beta| \|\exp(K\lambda\Omega)D^{\alpha-\beta}g\|_\infty.
		\end{equation*}
	\end{itemize}
\end{lemma}
\begin{prop}\label{inv swo}
	Let $\Omega$, $\Sigma$ be weight functions as in \eqref{directSum}; consider  $g\in \mathcal{S}^{\Sigma}_{\Omega}(\mathbb{R}^N)$, and a measurable function $F:\mathbb{R}^{2N}\to \mathbb{C}$ such that for each $\lambda>0$ there exists a constant $C_\lambda>0$ so that
	\begin{equation}\label{decay swo}
	|F(x,\xi)|\leq C_\lambda \exp(-\lambda(\Omega(x)+\Sigma(\xi))),
	\end{equation}
	for each $(x,\xi)\in \mathbb{R}^{2N}$. Then the integral
	\begin{equation}\label{2-6}
	f(t):=\int_{\mathbb{R}^{2N}} F(x,\xi) M_\xi T_x g(t) \, dx d\xi
	\end{equation}
	defines a function $f\in \mathcal{S}^{\Sigma}_{\Omega}(\mathbb{R}^N)$.
\end{prop}
\begin{pf}
	Observe that the integral in \eqref{2-6} is absolutely convergent, and we can differentiate under the integral sign. Fix $\lambda>0$ and $\alpha\in\mathbb{N}^N_0$. We get from Lemma $\ref{lemma gz}$
	\begin{align*}
		\|\exp(\lambda\Omega)D^\alpha f\|_\infty &\leq \int_{\mathbb{R}^{2N}} \|\exp(\lambda\Omega(t)) F(x,\xi) D^\alpha(M_\xi T_x g)(t)\|_{L^\infty(\mathbb{R}^N_t)} \, dx d\xi \\
		&\leq C\sum_{\beta \leq \alpha} \binom{\alpha}{\beta}\int_{\mathbb{R}^{2N}} \exp(K\lambda\Omega(x))|F(x,\xi)||\xi^\beta|\|\exp(K\lambda\Omega)D^{\alpha-\beta}g\|_\infty\, dx d\xi \\
&\leq C'_{\lambda,\alpha}\int_{\mathbb{R}^{2N}}|F(x,\xi)|P(x,\xi) \,dx d\xi,
	\end{align*}
	where $C'_{\lambda,\alpha}:=C\max\limits_{\beta\leq\alpha} \|\exp(K\lambda\Omega)D^{\alpha-\beta}g\|_\infty$ and
	\begin{equation*}
	P(x,\xi)=\sum_{\beta \leq \alpha} \binom{\alpha}{\beta} |\xi^\beta|\exp(K\lambda\Omega(x))=\exp(K\lambda\Omega(x))\prod_{j=1}^{N}(1+|\xi_j|)^{\alpha_j}.
	\end{equation*}
From \eqref{decay swo} we have that
\begin{equation*}
\int_{\mathbb{R}^{2N}} | F(x,\xi)| P(x,\xi)\,dxd\xi<\infty,
\end{equation*}
and so, for every $\lambda>0$ and $\alpha\in\mathbb{N}^N_0$,
\begin{equation}\label{finitenorm}
\|\exp(\lambda\Omega) D^\alpha f\|_\infty<\infty.
\end{equation}
Now we observe that, since for $x,\xi\in\mathbb{R}^N$, $\mathcal{F}(T_x g)=M_{-x}\hat{g}$ and $\mathcal{F}(M_\xi g)=T_\xi\hat{g}$,
	\begin{eqnarray*}
		\exp(\lambda\Sigma(t))D^\alpha \hat{f}(t)
&=&\int_{\mathbb{R}^{2N}} \exp(\lambda\Sigma(t)) F(x,\xi) D^\alpha_t\left( \exp(ix\xi) M_{-x}T_\xi\hat{g}\right)(t)\,dxd\xi,
\end{eqnarray*}
where we have used that $T_\xi M_{-x}=\exp(ix\xi) M_{-x} T_\xi$. Then, proceedings as before, we get
	\begin{equation*}
	\|\exp(\lambda\Sigma)D^\alpha \hat{f}\|_\infty  \leq D_{\lambda,\alpha}\int_{\mathbb{R}^{2N}}|F(x,\xi)|Q(x,\xi) \,dx d\xi,
	\end{equation*}
	with $D_{\lambda,\alpha}:=C\max\limits_{\beta \leq \alpha} \|\exp(K\lambda\Sigma)D^{\alpha-\beta}\hat{g}\|_\infty$ and
	\begin{equation*}
	Q(x,\xi)=\sum_{\beta \leq \alpha} \binom{\alpha}{\beta} |x^\beta|\exp(K\lambda\Sigma(\xi))=\exp(K\lambda\Sigma(\xi))\prod_{j=1}^{N}(1+|x_j|)^{\alpha_j}.
	\end{equation*}
Since $g\in\mathcal{S}_\Omega^\Sigma(\mathbb{R}^N)$ and $F$ satisfies \eqref{decay swo} we obtain
	\begin{eqnarray*}
		\| \exp(\lambda\Sigma)D^\alpha \hat{f}\|_\infty <\infty
	\end{eqnarray*}
for every $\lambda>0$ and $\alpha\in\mathbb{N}^N_0$ which, together with \eqref{finitenorm}, gives $f\in \mathcal{S}^{\Sigma}_{\Omega}(\mathbb{R}^N)$.\qed
\end{pf}
\begin{thm}\label{stime swo}
	Let $\Omega$, $\Sigma$ be weight functions as in \eqref{directSum} and consider  $g\in \mathcal{S}^{\Sigma}_{\Omega}(\mathbb{R}^N)$, $g\neq0$. Then for $f\in (\mathcal{S}^{\Sigma}_{\Omega})'(\mathbb{R}^N)$ the following conditions are equivalent:
	\begin{itemize}
		\item[(1)] $f\in \mathcal{S}^{\Sigma}_{\Omega}(\mathbb{R}^N)$;
		\item[(2)] $V_g f$ satisfies $(\ref{decay swo})$.
	\end{itemize}
\end{thm}
\begin{pf}
	$(1)\implies (2)$: fix $\lambda>0$. Then from Lemma \ref{prop swo}
	\begin{eqnarray*}
		\exp(2\lambda\Omega(x))|V_g f(x,\xi)|&\leq& \exp(2\lambda\Omega(x))\int_{\mathbb{R}^N} |\exp(-i t \xi) f(t)\overline{g(t-x)}|\, dt  \\
&\leq& C_{2\lambda} \| \exp(2K\lambda\Omega)f\|_\infty\| \exp(2K\lambda\Omega) \overline{Rg}\|_\infty <\infty,
	\end{eqnarray*}
	since $f,g \in \mathcal{S}^\Sigma_\Omega(\mathbb{R}^N)$. So for each $(x,\xi)\in \mathbb{R}^{2N}$ and for every $\lambda>0$,
	\begin{equation*}
	|V_g f(x,\xi)|\leq C'_\lambda\exp(-2\lambda\Omega(x)).
	\end{equation*}
	Analogously, using the foundamental identity of the STFT we get
	\begin{equation*}
		|V_gf(x,\xi)|=|(2\pi)^{-N}\exp(-ix\xi)V_{\hat{g}}\hat{f}(\xi,-x)| \leq  D_\lambda\exp(-2\lambda\Sigma(\xi)),
	\end{equation*}
where $D_\lambda=C_{2\lambda} \|\exp(2\lambda K\Sigma)\hat{f}\|_\infty \|\exp(2\lambda K\Sigma)\overline{R\hat{g}}\|_\infty$. Finally,
	\begin{equation*}
	|V_gf(x,\xi)|=\sqrt{|V_gf(x,\xi)|^2}\leq \sqrt{C_\lambda' D_\lambda} \exp(-\lambda\Omega(x))\exp(-\lambda\Sigma(\xi)).
	\end{equation*}
	$(2)\implies(1)$: from Proposition $\ref{inv swo}$, with $V_g f$ in place of $F$, and using the inversion
formula for the STFT (see for instance \cite{6}) we get $f\in \mathcal{S}^\Sigma_\Omega(\mathbb{R}^N)$.\qed
\end{pf}
\begin{oss}\label{in funx the}
	Observe that in the proof of Theorem $\ref{stime swo}$, when we prove that $(1)\implies(2)$, we only use the conditions:
	\begin{eqnarray*}
		\|\exp(\lambda\Omega)f\|_\infty <\infty, \quad \text{for each} \, \lambda>0,\\
		\| \exp(\lambda\Sigma)\hat{f}\|_\infty <\infty, \quad \text{for each} \, \lambda>0.
	\end{eqnarray*}
	Then if $f$ satisfies these conditions, $f\in \mathcal{S}^{\Sigma}_{\Omega}(\mathbb{R}^N)$.
\end{oss}

Now, for $\Omega$ and $\Sigma$ as in \eqref{directSum} we consider the Young conjugates $\varphi^*_{\omega_j}$ and $\varphi^*_{\sigma_j}$ of $\omega_j$ and $\sigma_j$ respectively, $j=1,\dots,N$, and we define the following functions for $y=(y_1,\dots,y_N)\in\mathbb{R}^N$, $y_j\geq 0$ for every $j=1,\dots,N$:
\begin{equation}\label{omegasigma*}
\Omega^*(y)=\varphi^*_{\omega_1}(y_1)+\dots+\varphi^*_{\omega_N}(y_N),\quad \Sigma^*(y)= \varphi^*_{\sigma_1}(y_1)+\dots+\varphi^*_{\sigma_N}(y_N).
\end{equation}
We have several equivalent conditions for $f$ to belong to $\mathcal{S}_\Omega^\Sigma(\mathbb{R}^N)$, summarized in the next result. Observe that the following theorem is proved, in the case of a single weight function, in \cite{8} (for $L^\infty$ norms) and in \cite{9} (for $L^p$-$L^q$ norms with $p,q<\infty$). Here we propose a unified proof for every $1\leq p,q\leq \infty$, in the case of non isotropic ultradifferentiable spaces.
\begin{thm}\label{semin ult swo}
	Let $\Omega$, $\Sigma$ be weight functions as in \eqref{directSum} and consider $f \in \mathcal{S}(\mathbb{R}^N)$; let moreover $1\leq p,q\leq\infty$. Then
the following conditions are equivalent:
	\begin{itemize}
\item[(1)] $f\in\mathcal{S}_\Omega^\Sigma(\mathbb{R}^N)$, i.e., it satisfies \eqref{cond 1 swo} and \eqref{cond 2 swo}.		
\item[(2)] $f$ satisfies the conditions:
		\begin{equation*}
		\|\exp(\lambda\Omega)D^\alpha f\|_p <\infty, \quad \text{for each} \, \lambda>0, \; \alpha\in\mathbb{N}^N_0,
		\end{equation*}
		\begin{equation*}
		\| \exp(\lambda\Sigma)D^\alpha \hat{f}\|_q <\infty, \quad \text{for each} \, \lambda>0,\; \alpha\in\mathbb{N}^N_0.
		\end{equation*}
		\item[(3)] $f$ satisfies the conditions:
		\begin{equation}\label{def oli 1}
		\|\exp(\lambda\Omega)f\|_p <\infty, \quad \text{for each} \, \lambda>0,
		\end{equation}
		\begin{equation}\label{def oli 2}
		\|\exp(\lambda\Sigma)\hat{f}\|_q <\infty, \quad \text{for each} \, \lambda>0.
		\end{equation}
		\item[(4)] $f$ satisfies the conditions:
		\begin{equation}\label{dd}
		\|\exp(\lambda\Omega)x^\alpha f\|_p <\infty, \quad \text{for each} \, \lambda>0, \; \alpha \in\mathbb{N}^N_0,
		\end{equation}
		\begin{equation}\label{parti}
		\|\exp(\lambda\Sigma)\xi^\alpha \hat{f}\|_q <\infty, \quad \text{for each} \, \lambda>0,\; \alpha \in\mathbb{N}^N_0.
		\end{equation}
		\item[(5)] $f$ satisfies the conditions:
		\begin{itemize}
			\item[(a)] For each $\lambda>0$ and each $\beta\in\mathbb{N}^N_0$ there exists $C_{\beta,\lambda}>0$ such that for each $\alpha\in\mathbb{N}^N_0$:
			\begin{equation}\label{5a} \left\|\exp\left(-\lambda\Sigma^*\left(\frac{\alpha}{\lambda}\right)\right) x^\beta D^\alpha f\right\|_p \leq C_{\beta,\lambda};
			\end{equation}
			\item[(b)] For each $\mu>0$ and each $\alpha\in\mathbb{N}^N_0$ there exists $C_{\alpha,\mu}>0$ such that for each $\beta\in\mathbb{N}^N_0$:
			\begin{equation}\label{5b} \left\|\exp\left(-\mu\Omega^*\left(\frac{\beta}{\mu}\right)\right)x^\beta D^\alpha f\right\|_q \leq C_{\alpha,\mu}.
			\end{equation}
		\end{itemize}
		\item[(6)] For each $\lambda,\mu>0$ there exists $C_{\mu,\lambda}>0$ such that for each $\alpha,\beta\in\mathbb{N}^N_0$:
		\begin{equation}\label{6} \left\|\exp\left(-\lambda\Sigma^*\left(\frac{\alpha}{\lambda}\right)-\mu\Omega^*\left(\frac{\beta}{\mu}\right)\right)x^\beta D^\alpha f\right\|_p \leq C_{\mu,\lambda}.
		\end{equation}
	\end{itemize}
\end{thm}
Before proving Theorem \ref{semin ult swo} we need two lemmas.
\begin{lemma}\label{infinity1}
Let $\Omega,\Sigma$ be weight functions as in \eqref{directSum} and $f\in\mathcal{S}(\mathbb{R}^N)$. If $f\in\mathcal{S}_\Omega^\Sigma(\mathbb{R}^N)$ then for each $\lambda,\mu>0$ there exists $C_{\mu,\lambda}>0$ such that for each $\alpha,\beta \in\mathbb{N}^N_0$
		\begin{equation}\label{3} \left\|\exp\left(-\lambda\Sigma^*\left(\frac{\alpha}{\lambda}\right)-\mu\Omega^*\left(\frac{\beta}{\mu}\right)\right)x^\beta D^\alpha f\right\|_\infty \leq C_{\mu,\lambda}.
		\end{equation}
\end{lemma}

\begin{pf}
We start by proving that for each $\mu>0$ and each $\alpha\in\mathbb{N}^N_0$ there exists $C_{\alpha,\mu}>0$ such that for each $\beta\in\mathbb{N}^N_0$:
			\begin{equation}\label{2bproof}
\left\|\exp\left(-\mu\Omega^*\left(\frac{\beta}{\mu}\right)\right)x^\beta D^\alpha f\right\|_\infty \leq C_{\alpha,\mu}.
			\end{equation}
From Proposition $\ref{wheight prop gamma'}$ (1) we have
	\begin{align*}
		 |x^\beta|\exp\left(-\mu\Omega^*\left(\frac{\beta}{\mu}\right)\right)\leq \exp(\mu\Omega(x)),
	\end{align*}
	and so
	\begin{align*}
		|x^\beta D^\alpha f(x)|\exp\left(-\mu\Omega^*\left(\frac{\beta}{\mu}\right)\right)\leq \exp(\mu\Omega(x))|D^\alpha f(x)|\leq \|\exp(\mu\Omega)D^\alpha f\|_\infty\leq C_{\alpha,\mu};
	\end{align*}
then \eqref{2bproof} holds. Now we want to prove that for each $\lambda>0$ and each $\beta\in\mathbb{N}^N_0$ there exists $C_{\beta,\lambda}>0$ such that for each $\alpha\in\mathbb{N}^N_0$:
			\begin{equation}\label{2aproof}
			 \left\|\exp\left(-\lambda\Sigma^*\left(\frac{\alpha}{\lambda}\right)\right)x^\beta D^\alpha f\right\|_\infty \leq C_{\beta,\lambda}.
			\end{equation}
We can write $x^\beta D^\alpha f(x)=\mathcal{F}^{-1}_{\xi\to x}\mathcal{F}_{t\to\xi}\left(t^\beta D^\alpha f(t)\right)$; then, using standard properties of the Fourier transform and Leibniz rule we have
	\begin{align*}
		|x^\beta D^\alpha f(x)|&\leq \sum_{\gamma\leq \alpha,\beta} \binom{\beta}{\gamma} (2\pi)^{-N} \int_{\mathbb{R}^{N}} |D_\xi^\gamma \xi^\alpha| |D_\xi^{\beta-\gamma}\hat{f}(\xi))| \, d\xi \\ &=
		\sum_{\gamma\leq \alpha,\beta} \frac{\beta!}{\gamma!(\beta-\gamma)!} (2\pi)^{-N} \int_{\mathbb{R}^{N}} \frac{\alpha!}{(\alpha-\gamma)!}|\xi^{\alpha-\gamma}||D_\xi^{\beta-\gamma}\hat{f}(\xi))| \, d\xi \\ &\leq 2^{|\alpha|}\sum_{\gamma\leq \alpha,\beta}\frac{(2\pi)^{-N}\beta!}{(\beta-\gamma)!} \int_{\mathbb{R}^{N}}  |\xi^{\alpha-\gamma}| |D_\xi^{\beta-\gamma}\hat{f}(\xi))| \exp(2\lambda'\Sigma(\xi)-2\lambda'\Sigma(\xi)) \, d\xi.
	\end{align*}

	By $(\ref{cond 2 swo})$ we have that for each $\gamma \leq \beta$
	\begin{equation*}
	\|\exp(2\lambda'\Sigma(\xi))D_\xi^{\beta-\gamma}\hat{f}\|_\infty\leq C'_{\beta,\lambda'}
	\end{equation*}
	for some $C'_{\beta,\lambda'}>0$. By Proposition $\ref{wheight prop gamma'}$ (1) we have that for each $\gamma \leq \alpha$
	\begin{equation*}
	|\xi^{\alpha-\gamma}|\exp(-\lambda'\Sigma(\xi)) \leq \exp\left(\lambda'\Sigma^*\left(\frac{\alpha}{\lambda'}\right)\right).
	\end{equation*}
	Therefore, we get
	\begin{align*}
		|x^\beta D^\alpha f(x)|\leq C''_{\beta,\lambda'}2^{|\alpha|}\exp\left(\lambda'\Sigma^*\left(\frac{\alpha}{\lambda'}\right)\right)\int_{\mathbb{R}^N} \exp(-\lambda'\Sigma(\xi))\, d\xi,
	\end{align*}
	for some $C''_{\beta,\lambda'}>0$. Using Proposition $\ref{wheight prop gamma'}$ (3), we obtain that for each $0\leq \lambda\leq \frac{\lambda'}{L^{[\log 2 +1]}}$ and for each $j=1,\dots,N$
	\begin{equation*}
	2^{|\alpha_j|} \exp\left(\lambda'\varphi^*_{\sigma_j}\left(\frac{|\alpha_j|}{\lambda'}\right)\right)\leq C_{\lambda'}\exp\left(\lambda\varphi^*_{\sigma_j}\left(\frac{|\alpha_j|}{\lambda}\right)\right),\end{equation*}
	for some $C_{\lambda'}>0$. Since $\lambda'$ is arbitrary we can choose also $\lambda$ in an arbitrary way; moreover, we can choose $\lambda'$ sufficiently large in such a way that $\int_{\mathbb{R}^N} \exp(-\lambda'\Sigma(\xi))\, d\xi <\infty$. We finally have that for each $\lambda>0$ there exists $C_{\beta,\lambda}>0$ such that
	\begin{equation*}
	|x^\beta D^\alpha f(x)|\leq C_{\beta,\lambda}\exp\left(\lambda\Sigma^*\left(\frac{\alpha}{\lambda}\right)\right),
	\end{equation*}
so we have proved \eqref{2aproof}. \\[0.1cm]
Now let us remark that, setting $\langle x\rangle=\sqrt{1+|x|^2}$ and $M=\left[\frac{N+1}{4}\right] +1$, we have
	\begin{align*}
		\|x^\beta D^\alpha f \|_2^2\leq C\|x^\beta \langle x \rangle^{2M} D^\alpha f \|_\infty^2.
	\end{align*}
	Since $\langle x \rangle^{2M}= \sum_{|\gamma| \leq M} \frac{M!}{\gamma!(M-|\gamma|)!} x^{2\gamma}$, we get
	\begin{eqnarray*}
		\|x^\beta D^\alpha f \|_2^2\leq C\sum_{|\gamma| \leq M} \frac{M!}{\gamma!(M-|\gamma|)!}\|x^{\beta+2\gamma} D^\alpha f \|_\infty^2.
	\end{eqnarray*}
	From condition $\eqref{2aproof}$, we then have that for each $\lambda>0$ and $\beta\in\mathbb{N}^N_0$ there exists $C_{\beta,\lambda}>0$ such that for each $\alpha\in\mathbb{N}^N_0$
	\begin{equation}\label{add1}
	\|x^\beta D^\alpha f \|_2 \leq C'_{\beta,\lambda}\exp\left(\lambda\Sigma^*\left(\frac{\alpha}{\lambda}\right)\right),
	\end{equation}
	where $C'_{\beta,\lambda}$ depends only on $\beta$, $\lambda$ and $M$. Analogously, from condition $\eqref{2bproof}$, by the convexity of the Young conjugate function, we get that for each $\mu'>0$ and $\alpha\in\mathbb{N}^N_0$ there exists $C_{\alpha,\mu'}>0$ such that for each $\beta\in\mathbb{N}^N_0$
	\begin{align*}
		\|x^\beta D^\alpha f \|_2^2 &\leq C\sum_{|\gamma| \leq M} \frac{M!}{|\gamma|!(M-|\gamma|)!}C^2_{\alpha,\mu'} \exp\left(2\mu'\Omega^*\left(\frac{\beta+2\gamma}{\mu'}\right)\right)\\ &\leq C'_{\alpha,\mu'}\exp\left(\mu'\Omega^*\left(\frac{2\beta}{\mu'}\right)\right),
 	\end{align*}
	where $C'_{\alpha,\mu'}$ depends only on $\alpha$, $\mu'$ and $M$. Then, setting $\mu=\frac{\mu'}{2}$, we get that for every $\alpha\in\mathbb{N}^N_0$ and $\mu>0$ there exists a constant $C''_{\alpha,\mu}>0$ such that for each $\beta\in\mathbb{N}^N_0$,
	\begin{equation}\label{add2}
	\|x^\beta D^\alpha f \|_2 \leq C''_{\alpha,\mu}\exp\left(\mu\Omega^*\left(\frac{\beta}{\mu}\right)\right).
	\end{equation}
	By integration by parts, Leibniz rule, H\"older's inequality, \eqref{add1} and \eqref{add2} we obtain
	\begin{align*}
		\|x^\beta D^\alpha f \|_2^2 &=\int_{\mathbb{R}^N} x^{2\beta} \overline{\partial^\alpha f(x)} \partial^\alpha f(x)\,dx \\
&\leq \sum_{\gamma\leq \alpha,2\beta} \binom{\alpha}{\gamma} \binom{2\beta}{\gamma} \gamma! \| \partial^{2\alpha-\gamma}f\|_2\| x^{2\beta-\gamma}f\|_2 \\
&\leq \sum_{\gamma\leq \alpha,2\beta} \binom{\alpha}{\gamma} \binom{2\beta}{\gamma} \gamma!C'_{0,\lambda}C''_{0,\mu}\exp\left(\lambda\Sigma^*\left(\frac{2\alpha-\gamma}{\lambda}\right)\right) \exp\left(\mu\Omega^*\left(\frac{2\beta-\gamma}{\mu}\right)\right).
	\end{align*}
	Now by Proposition $\ref{wheight prop gamma'}$ (4) and (3), and using the convexity of the Young conjugate we get
	\begin{align*}
		\|x^\beta D^\alpha f \|_2^2 &\leq \sum_{\gamma\leq \alpha,2\beta} \binom{\alpha}{\gamma} \binom{2\beta}{\gamma} C'_{0,\lambda}C''_{0,\mu}C_\lambda\exp\left(\lambda\Sigma^*\left(\frac{2\alpha}{\lambda}\right)\right) \exp\left(\mu\Omega^*\left(\frac{2\beta}{\mu}\right)\right)\\
&\leq C_{\lambda',\mu'}\exp\left(\lambda'\Sigma^*\left(\frac{2\alpha}{\lambda'}\right)\right) \exp\left(\mu'\Omega^*\left(\frac{2\beta}{\mu'}\right)\right).
	\end{align*}
	Then for each $\lambda,\mu>0$ there exists a constant $C'_{\lambda,\mu}>0$ such that
	\begin{equation}\label{add3}
	\|x^\beta D^\alpha f \|_2\leq C'_{\lambda,\mu}\exp\left(\lambda\Sigma^*\left(\frac{\alpha}{\lambda}\right)+\mu\Omega^*\left(\frac{\beta}{\mu}\right)\right).
	\end{equation}
	It is easy to show that the hypothesis $x^\beta D^\alpha f \in L^\infty(\mathbb{R}^N)$ for all $\alpha,\beta \in \mathbb{N}^N_0$ implies $x^\beta D^\alpha f \in H^s(\mathbb{R}^N)$ for all $\alpha$, $\beta\in\mathbb{N}^N_0$ and each $s>0$, where $H^s(\mathbb{R}^N)$ is the standard Sobolev space. Writing $\|^.\|_{H^s}$ for the Sobolev $H^s$ norm, by Sobolev inequality there exists $C'>0$ such that for $s>\frac{N}{2}$
	\begin{equation*}
	\|x^\beta D^\alpha f \|_\infty \leq C'  \|x^\beta D^\alpha f \|_{H^s}.
	\end{equation*}
	So in conclusion, fix an integer $s>\left[\frac{N}{2}\right]+1$. Then by \eqref{add3}
	\begin{align*}
		\|x^\beta D^\alpha f \|_\infty &\leq C' \sum_{|\gamma|\leq s}\Vert D^\gamma (x^\beta D^\alpha f)\Vert_2 \\
		&\leq C''\sum_{|\gamma|\leq s} \sum_{\delta\leq \gamma,\beta} \binom{\gamma}{\delta} \binom{\beta}{\delta} \delta!C'_{\lambda,\mu}\exp\left(\lambda\Sigma^*\left(\frac{\alpha+\gamma-\delta}{\lambda}\right)\right) \exp\left(\mu\Omega^*\left(\frac{\beta-\delta}{\mu}\right)\right).
	\end{align*}
	Proceeding as in the previous steps, using Proposition $\ref{wheight prop gamma'}$ (4) and the convexity of the Young conjugate function we get $(\ref{3})$. \qed
\end{pf}

\begin{lemma}\label{infinity2}
Let $\Omega,\Sigma$ be weight functions as in \eqref{directSum} and $f\in\mathcal{S}(\mathbb{R}^N)$. Suppose that $f$ satisfies the following conditions:
		\begin{itemize}
			\item[(a)] For each $\lambda>0$ and each $\beta\in\mathbb{N}^N_0$ there exists $C_{\beta,\lambda}>0$ such that for each $\alpha\in\mathbb{N}^N_0$:
			\begin{equation}\label{2a} \left\|\exp\left(-\lambda\Sigma^*\left(\frac{\alpha}{\lambda}\right)\right)x^\beta D^\alpha f\right\|_\infty \leq C_{\beta,\lambda};
			\end{equation}
			\item[(b)] For each $\mu>0$ and each $\alpha\in\mathbb{N}^N_0$ there exists $C_{\alpha,\mu}>0$ such that for each $\beta\in\mathbb{N}^N_0$:
			\begin{equation}\label{2b} \left\|\exp\left(-\mu\Omega^*\left(\frac{\beta}{\mu}\right)\right)x^\beta D^\alpha f\right\|_\infty \leq C_{\alpha,\mu}.
			\end{equation}
		\end{itemize}
Then $f\in\mathcal{S}_\Omega^\Sigma(\mathbb{R}^N)$.
\end{lemma}

\begin{pf}
We prove initially that $(\ref{cond 1 swo})$ holds. Let $\alpha\in\mathbb{N}^N_0$ and fix $\mu>0$. By \eqref{2b},
	\begin{align*}
		|D^\alpha f(x)| &\leq \left\|\exp\left(-\mu\Omega^*\left(\frac{\beta}{\mu}\right)\right)x^\beta D^\alpha f\right\|_\infty |x^{-\beta}|\exp\left(\mu\Omega^*\left(\frac{\beta}{\mu}\right)\right)  \\&\leq C_{\alpha,\mu}|x^{-\beta}|\exp\left(\mu\Omega^*\left(\frac{\beta}{\mu}\right)\right),
	\end{align*}
	for each $\beta\in\mathbb{N}^N_0$ and $x\neq 0$. Since $C_{\alpha,\mu}$ does not depend on $\beta$, we have that
	\begin{equation*}
	|D^\alpha f(x)|\leq C_{\alpha,\mu}\inf_{\beta\in \mathbb{N}_0^N} |x^{-\beta}|\exp\left(\mu\Omega^*\left(\frac{\beta}{\mu}\right)\right) .
	\end{equation*}
	Using Proposition $\ref{wheight prop gamma'}$ (2), we get that
    \begin{equation*}
        \|\exp(\mu\Omega)D^\alpha f\|_\infty\leq C'_{\alpha,\mu},
    \end{equation*}
	i.e., \eqref{cond 1 swo} is satisfied. Let us now prove $\eqref{cond 2 swo}$. Fix $\beta\in\mathbb{N}^N_0$; for $\xi\neq 0$ and $\alpha\in\mathbb{N}^N_0$ we have
	\begin{align*}
		|D^\beta \hat{f}(\xi)| &=\left\vert \int_{\mathbb{R}^N} x^\beta f(x)\exp(-ix\xi)\,dx\right\vert = \left\vert \int_{\mathbb{R}^N} \xi^{-\alpha} x^\beta f(x) D^\alpha_x\left(\exp(-ix\xi)\right)\,dx\right\vert \\
		&\leq \sum_{\gamma\leq \alpha,\beta} \frac{\beta!}{(\beta-\gamma)!} \binom{\alpha}{\gamma} \int_{\mathbb{R}^N} |x|^{|\beta-\gamma|}||D_x^{\alpha-\gamma}f(x)||\xi^{-\alpha}|\, dx
\\
&\leq  2^{|\alpha|}\sum_{\gamma\leq \alpha,\beta} \frac{\beta!}{(\beta-\gamma)!}  \int_{\mathbb{R}^N} \exp\left(-\lambda'\Sigma^*\left(\frac{\alpha-\gamma}{\lambda'}\right)\right) \exp\left(\lambda'\Sigma^*\left(\frac{\alpha-\gamma}{\lambda'}\right)\right) \times \\& \times \frac{\langle x\rangle^{|\beta-\gamma|+N+1}}{\langle x\rangle^{N+1}}|D_x^{\alpha-\gamma}f(x)||\xi^{-\alpha}|\, dx,
	\end{align*}
	where as usual $\langle x \rangle=\sqrt{1+|x|^2}$. By condition $(\ref{2a})$, we obtain that
	\begin{equation*}
	\langle x\rangle^{|\beta-\gamma|+N+1}|D_x^{\alpha-\gamma}f(x)| \exp\left(-\lambda'\Sigma^*\left(\frac{\alpha-\gamma}{\lambda'}\right)\right)\leq C_{\beta,\lambda'},
	\end{equation*}
	for some $C_{\beta,\lambda'}>0$. Moreover, applying Proposition $\ref{wheight prop gamma'}$ (3), we have that for each $0\leq \lambda\leq \frac{\lambda'}{L^{[\log 2 +1]}}$
	\begin{align*}
		|D^\beta \hat{f}(\xi)|\leq C'_{\beta,\lambda}\exp\left(\lambda\Sigma^*\left(\frac{\alpha}{\lambda}\right)\right)|\xi^{-\alpha}|,
	\end{align*}
	for some $C'_{\beta,\lambda'}>0$. Since this holds for each $\alpha\in\mathbb{N}^N_0$, we get from Proposition $\ref{wheight prop gamma'}$ (2)
	\begin{align*}
		|D^\beta \hat{f} (\xi)|&\leq C'_{\beta,\lambda}\inf_{\alpha\in \mathbb{N}_0^N} |\xi^{-\alpha}|\exp\left(\lambda\Sigma^*\left(\frac{\alpha}{\lambda}\right)\right)  \\&\leq  C_{\beta,\lambda}\exp\left(-\left(\lambda-\frac{1}{b}\right)\sigma_1(\xi)-\frac{a}{b}\right) \cdots\exp\left(-\left(\lambda-\frac{1}{b}\right)\sigma_N(\xi_N)-\frac{a}{b}\right),
	\end{align*}
	that implies
	\begin{equation*}
	\|\exp(\mu\Sigma)D^\beta \hat{f}\|_\infty\leq C^{''}_{\beta,\mu},
	\end{equation*}
	for some $C^{''}_{\beta,\mu}>0$, i.e., \eqref{cond 2 swo} is satisfied, and so $f\in\mathcal{S}_\Omega^\Sigma(\mathbb{R}^N)$.
\qed
\end{pf}

\begin{oss}
In view of Theorem \ref{semin ult swo} we have that, despite in Lemmas \ref{infinity1} and \ref{infinity2} we only prove an implication in one direction, they express in fact necessary and sufficient conditions. On the other hand, in order to prove Theorem \ref{semin ult swo} for every $1\leq p,q\leq\infty$, we only need, for $L^\infty$ norms, the results of Lemmas \ref{infinity1} and \ref{infinity2}.
\end{oss}

\begin{oss}\label{w in lp}
	Let $\omega$ be a weight function as in Definition \ref{weight} and consider $1\leq p <\infty$. Then $\exp(-\lambda\omega)\in L^p(\mathbb{R}^N)$ if and only if $\lambda> \frac{N}{bp}$, where $b$ is the constant appearing in condition $(\gamma)$ of Definition \ref{weight}.
\end{oss}

\begin{pf}[of Theorem \ref{semin ult swo}]
	$(3) \implies (1)$:  we want to use Theorem $\ref{stime swo}$. In particular, we fix $g\in \mathcal{S}^{\Sigma}_{\Omega}(\mathbb{R}^N)$, $g\neq 0$ and we prove that $(\ref{decay swo})$ holds for $V_gf(x,\xi)$. Consider $r> \frac{N}{bp'}$ where $p'$ is the conjugate exponent of $p$, (if $p=1$ we choose any $r>0$). Using H\"older inequality we have
	\begin{align*}
		|\exp(\lambda\Omega(x))V_g f(x,\xi)|\leq \exp(\lambda NK)\|\exp((\lambda K+r) \Omega)g\|_\infty \|\exp(\lambda K \Omega) f \|_p \|\exp(-r\Omega)\|_{p'} <\infty,
	\end{align*}
	thanks to Remark $\ref{w in lp}$ and since $f$ satisfies $\eqref{def oli 1}$, $g$ satisfies $\eqref{cond 1 swo}$. Analogously, using the foundamental identity of the STFT and $q$ instead of $p$, we get, for $r'> \frac{N}{bq'}$
	\begin{align*}
		&|\exp(\lambda\Sigma(\xi))V_g f(x,\xi)|\leq \\
&\quad\leq (2\pi)^{-N}\exp(\lambda NK)\|\exp((\lambda K+r') \Sigma)\hat{g}\|_\infty
		\|\exp(\lambda K \Sigma) \hat{f} \|_q \|\exp(-r'\Sigma)\|_{q'} <\infty,
	\end{align*}
	thanks to Remark $\ref{w in lp}$ and since $f$ satisfies $\eqref{def oli 2}$, $g$ satisfies $\eqref{cond 2 swo}$.
	\\[0.1cm]
	Hence we have that for each $\lambda>0$
	\begin{equation*}
	|V_gf(x,\xi)|=\sqrt{|V_gf(x,\xi)|^2}\leq C_\lambda \exp(-\lambda\Omega(x))\exp(-\lambda\Sigma(\xi)).
	\end{equation*}
	for some $C_\lambda>0$. From Theorem \ref{stime swo}, we get $f\in\mathcal{S}_\Omega^\Sigma(\mathbb{R}^N)$.
	\\[0.2cm]
	$(2)\implies (3)$: it is trivial, taking $\alpha=0$.
	\\[0.2cm]
	$(1)\implies (2)$: we suppose that $p,q<\infty$ (otherwise the corresponding implication is trivial). Fix $\lambda>0$, $\alpha\in\mathbb{N}^N_0$ and consider $r> \frac{N}{bp}$. We have
	\begin{align*}
		\|\exp(\lambda\Omega)D^\alpha f\|_p^p \leq \|\exp((\lambda+r)\Omega)D^\alpha f\|_\infty^p \|\exp(-r\Omega)\|_p^p <\infty.
	\end{align*}
	In the same way, setting $r'> \frac{N}{bq}$
	\begin{align*}
		\|\exp(\lambda\Sigma)D^\alpha\hat{f}\|_q^q\leq \|\exp((\lambda+r')\Sigma)D^\alpha\hat{f}\|_\infty^q \|\exp(-r'\Sigma)\|_q^q <\infty.
	\end{align*}
	\\[0.2cm]
	$(4)\implies (3)$: it is trivial, taking $\alpha=0$.
	\\[0.2cm]
	$(3)\implies (4)$: fix $\alpha\in\mathbb{N}^N_0$ and $\lambda>0$. From condition $(\gamma)$ we get
	\begin{eqnarray*}
		\exp(\lambda\Omega(x))|x^\alpha| \leq C_\alpha \exp(\lambda'\Omega(x)).
	\end{eqnarray*}
	Then
	\begin{equation*}
	\|\exp(\lambda\Omega)x^\alpha f\|_p \leq C_\alpha\|\exp(\lambda'\Omega) f\|_p <\infty,
	\end{equation*}
	since $f$ satisfies $\eqref{def oli 1}$. Analogously
	\begin{equation*}
	\|\exp(\lambda\Sigma)\xi^\alpha \hat{f}\|_q \leq C_\alpha\|\exp(\lambda'\Sigma) \hat{f}\|_q <\infty,
	\end{equation*}
	since $f$ satisfies $\eqref{def oli 2}$.
	\\[0.2cm]
	$(5)\implies (1)$: by Remark \ref{s} it is sufficient to show that $\hat{f} \in \mathcal{S}^{\Omega}_{\Sigma}(\mathbb{R}^N)$. In order to do this, we shall prove that $\hat{f}$ satisfies the hypotheses of Lemma \ref{infinity2} with $\Omega$ and $\Sigma$ interchanged. Fix $M>\left[{\frac{N}{2p'}}\right]+1$. For each $\alpha, \beta\in\mathbb{N}^N_0$ and each $\xi \in \mathbb{R}^N$ we have
	\begin{align*}
		|\xi^\beta D^\alpha \hat{f}(\xi)| \leq \|\langle x \rangle^{-2M}\|_{p'} \|\langle x \rangle^{2M}D_x^\beta(x^\alpha f)\|_p   \leq C_M \sum_{\gamma\leq \alpha,\beta} \binom{\beta}{\gamma} \binom{\alpha}{\gamma}\gamma!\|\langle x \rangle^{2M}x^{\alpha-\gamma}D^{\beta-\gamma} f\|_p
	\end{align*}
	for some $C_M>0$. Since $\langle x \rangle^{2M}= \sum_{|\delta| \leq M} \frac{M!}{\delta!(M-|\delta|)!} x^{2\delta}$, substituting in the previous estimate and using \eqref{5a}, we obtain that for each $\lambda>0$ and each $\alpha,\gamma,\delta\in\mathbb{N}^N_0$ there exists $C_{\alpha-\gamma+2\delta,\lambda}>0$ such that
	\begin{equation*}
		|\xi^\beta D^\alpha \hat{f}(\xi)|\leq C_M \sum_{\gamma\leq \alpha,\beta}\sum_{|\delta| \leq M} \frac{M!}{\delta!(M-|\delta|)!} \binom{\beta}{\gamma} \binom{\alpha}{\gamma}\gamma!C_{\alpha-\gamma+2\delta,\lambda} \exp\left(\lambda\Sigma^*\left(\frac{\beta-\gamma}{\lambda}\right)\right).
	\end{equation*}
	From the convexity of the Young conjugate function and by Proposition $\ref{wheight prop gamma'}$ (4), we see that
	\begin{align*}
		|\xi^\beta D^\alpha \hat{f}(\xi)|&\leq C_M \sum_{\gamma\leq \alpha,\beta}\sum_{|\delta| \leq M} \frac{M!}{\delta!(M-|\delta|)!} \binom{\beta}{\gamma} \binom{\alpha}{\gamma} C_\lambda\exp\left(\lambda\Sigma^*\left(\frac{\gamma}{\lambda}\right)\right) \times\\&\times C_{\alpha-\gamma+2\delta,\lambda}\exp\left(\lambda\Sigma^*\left(\frac{\beta-\gamma}{\lambda}\right)\right).
	\end{align*}
	From the convexity of $\varphi^*_{\sigma_j}$ and by Proposition $\ref{wheight prop gamma'}$ (3), we finally get
	\begin{align*}
		|\xi^\beta D^\alpha \hat{f}(\xi)|&\leq C'_{\alpha,\lambda'} \exp\left(\lambda'\Sigma^*\left(\frac{\beta}{\lambda'}\right)\right),
	\end{align*}
	for each $\lambda'>0$, $\alpha\in\mathbb{N}^N_0$, and for some $C'_{\alpha,\lambda'}>0$. In particular, $\hat{f}$ satisfies \eqref{2b} with $\Sigma^*$ in place of $\Omega^*$. In the same way, from \eqref{5b} we get that $\hat{f}$ satisfies \eqref{2a} with $\Omega^*$ in place of $\Sigma^*$, and so by Lemma \ref{infinity2} and Remark \ref{s} we have the claim.
	\\[0.2cm]
	$(1.)\implies (6.)$: from Lemma \ref{infinity1} we have that  $f$ satisfies $\eqref{3}$. Fix $\lambda,\mu>0$ and consider $M>\left[\frac{N}{2p}\right]+1$ and $\alpha,\beta\in\mathbb{N}^N_0$. Observe that
\begin{equation*}
\Vert x^\beta D^\alpha f\Vert_p\leq \Vert \langle x\rangle^{-2M}\Vert_p \Vert \langle x\rangle^{2M} x^\beta D^\alpha f\Vert_\infty;
\end{equation*}
then, applying the same technique as above, we have that $f$ satisfies \eqref{6}
	\\[0.2cm]
	$(6.)\implies (5.)$: trivial.
\qed
\end{pf}

\begin{prop}
	Let $\Omega$, $\Sigma$ be weight functions as in \eqref{directSum}. Then $\mathcal{S}^{\Sigma}_{\Omega}(\mathbb{R}^N)$ is closed under differentiation $D^\alpha$ and  multiplication by $x^\alpha$, with $\alpha\in\mathbb{N}^N_0$.
\end{prop}
\begin{pf}
	Let $\alpha\in\mathbb{N}^N_0$ and fix $f \in \mathcal{S}^\Sigma_\Omega(\mathbb{R}^N)$. Since $\mathcal{F}(D^{\alpha}f)(\xi)= \xi^\alpha \hat{f}(\xi)$ and $D^\beta (D^\alpha f)(x)=D^{\alpha+\beta}f(x)$, thanks to Theorem $\ref{semin ult swo}$ we get that $D^\alpha f \in \mathcal{S}^\Sigma_\Omega(\mathbb{R}^N)$. Analogously, $\mathcal{F}(x^{\alpha}f)(\xi)=(-1)^{|\alpha|} D^\alpha \hat{f}(\xi)$ and the same theorem implies that $x^\alpha f \in \mathcal{S}^\Sigma_\Omega(\mathbb{R}^N)$. \qed
\end{pf}

\section{Time-frequency representations and non isotropic ultradifferentiable classes}\label{h}

In this section we analyze the action of some transformations, namely the Short-time Fourier transform and a transform of Wigner type, on the space $\mathcal{S}_\Omega^\Sigma$. Concerning the Short-time Fourier transform we have the following remark.
\begin{oss}\label{vgf}
	Let $\Omega$ and $\Sigma$ be weight functions in $\mathbb{R}^N$ as in \eqref{directSum}. By Theorem $\ref{stime swo}$, we know that for a non-zero window $g\in \mathcal{S}^{\Sigma}_{\Omega}(\mathbb{R}^N)$ and for $f\in \mathcal{S}^{\Sigma}_{\Omega}(\mathbb{R}^N)$ we have
	\begin{equation}\label{add5}
	|V_g f(x,\xi)|\leq C_\lambda \exp(-\lambda(\Omega(x)+\Sigma(\xi)))
	\end{equation}
	for each $\lambda>0$. Now we observe that since $ g \in \mathcal{S}^{\Sigma}_{\Omega}(\mathbb{R}^N)$ then, from Remark \ref{s}, $\hat{g} \in \mathcal{S}^{\Omega}_{\Sigma}(\mathbb{R}^N)$. Hence for each $\lambda>0$ we have that
	\begin{equation}\label{add4}
	|\mathcal{F}V_g f (y,\eta)|=|(2\pi)^Nf(-\eta)\hat{g}(y)|\leq
	C'_\lambda\exp(-\lambda \Omega(\eta))C''_\lambda\exp(-\lambda \Sigma(y)).
	\end{equation}
Then by Theorem \ref{semin ult swo} we have that $V_g f\in\mathcal{S}_{\Omega\oplus\Sigma}^{\Sigma\oplus\Omega}(\mathbb{R}^{2N})$.

\end{oss}

Observe that the Short-time Fourier transform is defined as a partial Fourier transform of a function in $\mathbb{R}^{2N}$. We now analyze the action of partial Fourier transforms on ultradifferentiable spaces. Given $F\in \mathcal{S}(\mathbb{R}^{2N})$, we denote by $\mathcal{F}_1 F$ (resp. $\mathcal{F}_2 F$) the partial Fourier transform with respect to the first (resp. second) $N$-variables; explicitly,
	\begin{align*}
		&\mathcal{F}_1 F(y,\xi):=\int_{\mathbb{R}^N}\exp(-ixy)F(x,\xi)\, dx, \\&
		\mathcal{F}_2 F(x,\eta):=\int_{\mathbb{R}^N}\exp(-i\xi\eta)F(x,\xi)\, d\xi.
	\end{align*}

Observe that $\mathcal{F}=\mathcal{F}_2\mathcal{F}_1=\mathcal{F}_1\mathcal{F}_2$ and hence $\mathcal{F}^{-1}=\mathcal{F}^{-1}_1\mathcal{F}^{-1}_2=\mathcal{F}^{-1}_2\mathcal{F}^{-1}_1$.
\\[0.2cm]
Now we consider a collection of weight functions $\omega_{i,j}$, $\sigma_{i,j}$ for $i\in\{1,2\}$ and $j=1,\dots,N$; we write

\begin{equation}\label{Omega}
\Omega(x,y)=\omega_{1,1}(x_1)+\dots+\omega_{1,N}(x_N)+\omega_{2,1}(y_1)+\dots+\omega_{2,N}(y_N),
\end{equation}
\begin{equation}\label{Sigma}
\Sigma(\xi,\eta)=\sigma_{1,1}(\xi_1)+\dots+\sigma_{1,N}(\xi_N)+\sigma_{2,1}(\eta_1)+\dots+\sigma_{2,N}(\eta_N).
\end{equation}
\begin{equation}\label{Omega1}
\Omega_1(x,y)=\omega_{1,1}(x_1)+\dots+\omega_{1,N}(x_N)+\sigma_{2,1}(y_1)+\dots+\sigma_{2,N}(y_N),
\end{equation}
\begin{equation}\label{Sigma1}
\Sigma_1(\xi,\eta)=\sigma_{1,1}(\xi_1)+\dots+\sigma_{1,N}(\xi_N)+\omega_{2,1}(\eta_1)+\dots+\omega_{2,N}(\eta_N).
\end{equation}
Observe that we use the same notation $\Omega$ and $\Sigma$ for weights of the form \eqref{directSum} in $\mathbb{R}^N$ and of the form \eqref{Omega}-\eqref{Sigma} in $\mathbb{R}^{2N}$, since they are the same object. In $\mathbb{R}^{2N}$ we shall need in the following the associated weights $\Omega_1$, $\Sigma_1$; it will always be clear what we mean when we write $\Omega$, $\Sigma$.

\begin{thm}\label{pft swo}
	Let $\Omega$, $\Sigma$ be weight functions in $\mathbb{R}^{2N}$ as in \eqref{Omega} and \eqref{Sigma}. We have the following continuous maps.
	\begin{itemize}
		\item[(1)] $\mathcal{F}_2:\mathcal{S}^{\Sigma}_{\Omega}(\mathbb{R}^{2N})\to\mathcal{S}_{\Omega_1}^{\Sigma_1}(\mathbb{R}^{2N})$;
\item[(2)]
$\mathcal{F}_1:\mathcal{S}^{\Sigma}_{\Omega}(\mathbb{R}^{2N})\to\mathcal{S}_{\Sigma_1}^{\Omega_1}(\mathbb{R}^{2N})$.
   \end{itemize}
\end{thm}
\begin{pf}
The proofs of the two points are very similar; we show in detail point (1).
\\[0.1cm]
 Consider $F \in \mathcal{S}^{\Sigma}_{\Omega}(\mathbb{R}^{2N})$. Then $F$ satisfies conditions $(\ref{def oli 1})$ and $(\ref{def oli 2})$. Now fix $\lambda>0$ sufficiently large (of course it is sufficient to prove the corresponding estimates on $\mathcal{F}_2 f$ for $\lambda$ large). We have
	\begin{align*}
		|\mathcal{F}_2 F(x,\eta)| &\leq \int_{\mathbb{R}^N} |F(x,\xi)|\, d\xi  \leq \int_{\mathbb{R}^N} C_{2\lambda}\exp(-2\lambda(\Omega(x,\xi)))\, d\xi \\
&= C'_{2\lambda} \exp(-2\lambda(\omega_{1,1}(x_1)+\dots+\omega_{1,N}(x_N))).
	\end{align*}
	Moreover
	\begin{align*}
		|\mathcal{F}_2 F(x,\eta)|&=|\mathcal{F}_1^{-1}\mathcal{F}F(x,\eta)| \leq (2\pi)^{-N} \int_{\mathbb{R}^N} |\mathcal{F}F(t,\eta)|\, dt \\
&\leq (2\pi)^{-N} \int_{\mathbb{R}^N} D_{2\lambda}\exp(-2\lambda(\Sigma(t,\eta))) \, dt  =D'_{2\lambda} \exp(-2\lambda(\sigma_{2,1}(\eta_1)+\dots+\sigma_{2,N}(\eta_N))).
	\end{align*}
	Therefore, we have
	\begin{equation*}
	|\mathcal{F}_2 F(x,\eta)|=\sqrt{|\mathcal{F}_2 F(x,\eta)|^2}\leq \sqrt{C'_{2\lambda}D'_{2\lambda}} \exp(-\lambda\Omega_1 (x,\eta)).
	\end{equation*}
	Now similarly, we obtain
	\begin{eqnarray*}
		|\mathcal{F}\mathcal{F}_2F(y,\eta)|&=&|(2\pi)^N \mathcal{F}_1 F(y,-\eta)|
\leq (2\pi)^N \int_{\mathbb{R}^N} |F(x,-\eta)| \, dx \\
&\leq& (2\pi)^N \int_{\mathbb{R}^N} C_{2\lambda}\exp(-2\lambda(\Omega(x,\eta))) \, dx \\
&=& C'_{2\lambda} \exp(-2\lambda(\omega_{2,1}(\eta_1)+\dots+\omega_{2,N}(\eta_N))),
	\end{eqnarray*}
	and
	\begin{eqnarray*}
		|\mathcal{F}\mathcal{F}_2F(y,\eta)|&=&|\mathcal{F}_2 \mathcal{F} F(y,\eta)|
\leq  \int_{\mathbb{R}^N} |\mathcal{F}F(y,\xi)| \, d\xi \\
&\leq& \int_{\mathbb{R}^N} D_{2\lambda}\exp(-2\lambda(\Sigma(y,\xi))) \, d\xi   =D'_{2\lambda} \exp(-2\lambda(\sigma_{1,1}(y_1)+\dots+\sigma_{1,N}(y_N))).
	\end{eqnarray*}
	Hence
	\begin{equation*}
	|\mathcal{F} \mathcal{F}_2 F(y,\eta)|=\sqrt{|\mathcal{F}\mathcal{F}_2 F(y,\eta)|^2}\leq \sqrt{C'_{2\lambda}D'_{2\lambda}} \exp(-\lambda\Sigma_1(y,\eta)).
	\end{equation*}
	Then by Theorem \ref{semin ult swo} we get $\mathcal{F}_2 F\in\mathcal{S}_{\Omega_1}^{\Sigma_1}(\mathbb{R}^{2N})$.
\qed
\end{pf}
\begin{oss}
	We have already shown in Remark \ref{vgf}, that if $f,g \in \mathcal{S}^{\Sigma}_{\Omega}(\mathbb{R}^N)$, for $\Omega$ and $\Sigma$ as in \eqref{directSum}, then $V_g f \in \mathcal{S}_{\Omega\oplus\Sigma}^{\Sigma\oplus\Omega}(\mathbb{R}^{2N})$. We can re-obtain this result by means of Theorem \ref{pft swo}, since $V_g f(x,\xi)= (\mathcal{F}_2 F)(x,\xi)$, where $F(x,t)=f(t)\overline{g(t-x)}$ belongs to $\mathcal{S}^{\Sigma\oplus\Sigma}_{\Omega\oplus\Omega}(\mathbb{R}^{2N})$.
\end{oss}
By standard duality arguments we have the following result.
\begin{thm}\label{pft swo dis}
	Let $\Omega$, $\Sigma$ be weight functions as in \eqref{Omega}, \eqref{Sigma}. Then
	\begin{itemize}
		\item[(1)] $\mathcal{F}_2:(\mathcal{S}^{\Sigma}_{\Omega})'(\mathbb{R}^{2N})\to\left(\mathcal{S}_{\Omega_1}^{\Sigma_1}\right)'(\mathbb{R}^{2N})$;
	\item[(2)]
$\mathcal{F}_1:(\mathcal{S}^{\Sigma}_{\Omega})'(\mathbb{R}^{2N})\to\left(\mathcal{S}_{\Sigma_1}^{\Omega_1}\right)'(\mathbb{R}^{2N})$.

	\end{itemize}
\end{thm}
Now we introduce a Wigner-like transform and study the corresponding mapping properties in weighted ultradifferentiable spaces, as a preparation for applications to PDEs that we develop in the last part of the paper.
\begin{defn}
	Given $u\in \mathcal{S}(\mathbb{R}^{2N})$, we define $Wig[u]: \mathbb{R}^{2N}\to \mathbb{C}$ as
	\begin{equation}\label{cc}
	Wig[u](x,\xi)=\int_{\mathbb{R}^N}\exp(-it\xi)u\left(x+\frac{t}{2},x-\frac{t}{2}\right)\, dt.
	\end{equation}
\end{defn}
\begin{oss}\label{symmcoord}
	Let us consider the symmetric coordinate change $\mathfrak{T}_z$ acting on a function $F$ on $\mathbb{R}^{2N}$ as
	\begin{equation*}
	\mathfrak{T}_z F(x,\xi)= F\left(x+\frac{\xi}{2},x-\frac{\xi}{2}\right).
	\end{equation*}
	$\mathfrak{T}_z$ is a linear and bounded operator on $L^2(\mathbb{R}^{2N})$. Moreover it is invertible with inverse given by
	\begin{equation*}
	\mathfrak{T}_z^{-1} F(x,\xi)= F\left(\frac{x+\xi}{2},x-\xi\right).
	\end{equation*}
We observe that for $u\in\mathcal{S}(\mathbb{R}^{2N})$, $Wig[u]=\mathcal{F}_2\mathfrak{T}_z u$, and since both $\mathcal{F}_2$ and $\mathfrak{T}_z$ can be extended in a standard way to  (ultra)distributions, we can consider $Wig[u]$ acting on ultradistributions. Moreover, $Wig$ is invertible with inverse $Wig^{-1}= \mathfrak{T}_z ^{-1}\mathcal{F}_2^{-1}$.
\end{oss}
\begin{prop}\label{sim}
	Let $\Omega$, $\Sigma$ be weight functions as in \eqref{Omega}, \eqref{Sigma}. Then the following properties hold:
	\begin{itemize}
		\item[(1)]
$Wig: \mathcal{S}(\mathbb{R}^{2N})\to \mathcal{S}(\mathbb{R}^{2N})$;
		\item[(2)]
$Wig: \mathcal{S}'(\mathbb{R}^{2N})\to \mathcal{S}'(\mathbb{R}^{2N})$;
		\item[(3)]
$Wig: \mathcal{S}^{\Sigma}_{\Omega}(\mathbb{R}^{2N})\to\mathcal{S}_{\Omega_1}^{\Sigma_1}(\mathbb{R}^{2N})$;
		\item[(4)]
$Wig: (\mathcal{S}^{\Sigma}_{\Omega})'(\mathbb{R}^{2N})\to \left(\mathcal{S}_{\Omega_1}^{\Sigma_1}\right)'(\mathbb{R}^{2N})$.
	\end{itemize}
\end{prop}
\begin{pf}
It follows from Remark \ref{symmcoord}, Theorem \ref{pft swo} and Theorem \ref{pft swo dis}. \qed
\end{pf}

Similarly we have the following result.
\begin{prop}\label{flo}
	Let $\Omega$, $\Sigma$ be weight functions as in \eqref{Omega},\eqref{Sigma}. Then the following properties hold:
		\begin{itemize}
		\item[(1)] $Wig^{-1}: \mathcal{S}(\mathbb{R}^{2N})\to \mathcal{S}(\mathbb{R}^{2N})$;
		\item[(2)] $Wig^{-1}: \mathcal{S}'(\mathbb{R}^{2N})\to \mathcal{S}'(\mathbb{R}^{2N})$;
		\item[(3)] $Wig^{-1}: \mathcal{S}^{\Sigma}_{\Omega}(\mathbb{R}^{2N})\to\mathcal{S}_{\Omega_1}^{\Sigma_1}(\mathbb{R}^{2N})$;
		\item[(4)] $Wig^{-1}: (\mathcal{S}^{\Sigma}_{\Omega})'(\mathbb{R}^{2N})\to \left(\mathcal{S}_{\Omega_1}^{\Sigma_1}\right)'(\mathbb{R}^{2N})$.
	\end{itemize}
\end{prop}

\section{Regularity of PDE in weighted ultradifferentiable spaces}

In this section we give an application of the results that we have proved to the theory of PDEs. In particular we analyze the problem of regularity of solutions of partial differential equations with polynomial coefficients, and see how representation of time-frequency type can be profitably used in this field. Results in this direction have been studied in \cite{12} and \cite{8}, in the classical Schwartz space and in isotopic ultradifferentiable classes, in dimension $1$; here we provide a general framework for anisotropic spaces in arbitrary dimension $N$.

In order to state or results we need some notations. Let $\mathbf{R}:\mathcal{B}\to\mathcal{B}^L$, $L\in\mathbb{N}$, be a vector of $L$ operators acting on a space $\mathcal{B}$; then there exist $R_j:\mathcal{B}\to \mathcal{B}$, $j=1,\dots,L$, and for every $w\in \mathcal{B}$ we have
$$
\mathbf{R}w=(R_1 w,\dots,R_L w).
$$
In the following $\mathcal{B}$ will be either $C^\infty(\mathbb{R}^M)$, $\mathcal{S}(\mathbb{R}^M)$, or $\mathcal{S}_\Omega^\Sigma(\mathbb{R}^M)$. Then for every $\gamma\in\mathbb{N}^L_0$, following the multi-index notation we define the operator $\mathbf{R}^\gamma: \mathcal{B}\to \mathcal{B}$ as the composition
$$
\mathbf{R}^\gamma  = R_1^{\gamma_1}\cdots R_L^{\gamma_L}.
$$
Observe that for every $j=1,\dots,L$ we have $R_j=\mathbf{R}^{e_j}$, where $e_j$, $j=1,\dots,L$, is the $j$-th vector of the canonical basis of $\mathbb{R}^L$ (i.e., the vector having $1$ in the $j$-th position and $0$ elsewhere), so we can write
$$
\mathbf{R}=(\mathbf{R}^{e_1},\dots,\mathbf{R}^{e_L}).
$$
If for every $i,j=1,\dots,L$ the commutation relation
\begin{equation}\label{commutGeneral}
R_iR_j=R_jR_i
\end{equation}
holds, then for every $\gamma,\mu\in\mathbb{N}^L_0$ we have
\begin{equation}\label{AddSumOper}
\mathbf{R}^\gamma \mathbf{R}^\mu=\mathbf{R}^{\gamma+\mu};
\end{equation}
Now let $\mathbf{T}=(T_1,\dots,T_L)$ be another vector of $L$ operators acting on $\mathcal{B}$, and $a,b\in\mathbb{C}$; then for every $\gamma\in\mathbb{N}^L_0$ the operator $(a\mathbf{R}+b\mathbf{T})^\gamma$ is the composition
\begin{equation}\label{not}
(a\mathbf{R}+b\mathbf{T})^\gamma = (aR_1+bT_1)^{\gamma_1}\cdots (aR_L+bT_L)^{\gamma_L}.
\end{equation}
\begin{defn}\label{defMD}
	Let us consider a function $F\in C^\infty (\mathbb{R}^{2N})$, $F=F(x,\xi)$ for $x,\xi\in\mathbb{R}^N$. We define the multiplication and differentiation with respect to the first and second $N$ variables as the vectors of operators acting on $F(x,\xi)$ as
	\begin{eqnarray*}
	\mathbf{M_f} F(x,\xi)=(x_1 F(x,\xi),\dots,x_N F(x,\xi)),\\
	\mathbf{M_s} F(x,\xi)=(\xi_1 F(x,\xi),\dots,\xi_N F(x,\xi)),\\
	\mathbf{D_f} F(x,\xi)= (D_{x_1} F(x,\xi),\dots,D_{x_N}F(x,\xi)),\\
	\mathbf{D_s} F(x,\xi)= (D_{\xi_1} F(x,\xi),\dots,D_{\xi_N}F(x,\xi)).
	\end{eqnarray*}
\end{defn}

\begin{oss}
For a function $F\in C^\infty (\mathbb{R}^{2N})$, $F=F(x,\xi)$ for $x,\xi\in\mathbb{R}^N$, and a multi-index $\alpha\in\mathbb{N}^N_0$ we have
	\begin{eqnarray*}
	&\mathbf{M}_{\mathbf{f}}^\alpha F(x,\xi)=x^\alpha F(x,\xi),\quad
	&\mathbf{M}_{\mathbf{s}}^\alpha F(x,\xi)=\xi^\alpha F(x,\xi),\\
	&\mathbf{D}_{\mathbf{f}}^\alpha F(x,\xi)= D_x^\alpha F(x,\xi),\quad
	&\mathbf{D}_{\mathbf{s}}^\alpha F(x,\xi)= D_\xi^\alpha F(x,\xi).
	\end{eqnarray*}
\end{oss}
Now consider an operator with polynomial coefficients
\begin{equation}\label{ante}
P(x,y,D_x,D_y)=\sum_{|\alpha+\beta+\gamma+\mu|\leq m} c_{\alpha\beta\gamma\mu}x^\alpha y^\beta D_x^\gamma D_y^\mu,
\end{equation}
with $x,y\in\mathbb{R}^N$, $c_{\alpha\beta\gamma\mu}\in\mathbb{C}$, and $m\in\mathbb{N}$. If, for $\mathbf{j}=\mathbf{1},\dots,\mathbf{4}$, the operators $\mathbf{A^{j}}$ and $\mathbf{B^{j}}$ are any of $\mathbf{M_f}, \mathbf{M_s}, \mathbf{D_f}, \mathbf{D_s}$ and $a_j,b_j\in\mathbb{R}$, $j=1,\dots,4$, we denote
\begin{equation}\label{not1}
\begin{split}
&P\left(a_1\mathbf{A^{1}}+b_1\mathbf{B^{1}},a_2\mathbf{A^{2}}+b_2\mathbf{B^{2}}, a_3\mathbf{A^{3}}+b_3\mathbf{B^{3}},a_4\mathbf{A^{4}}+b_4\mathbf{B^{4}}\right)= \\
&\quad =\sum_{|\alpha+\beta+\gamma+\mu|\leq m} c_{\alpha\beta\gamma\mu} (a_1\mathbf{A^{1}}+b_1\mathbf{B^{1}})^\alpha (a_2\mathbf{A^{2}}+b_2\mathbf{B^{2}})^\beta (a_3\mathbf{A^{3}}+b_3\mathbf{B^{3}})^\gamma (a_4\mathbf{A^{4}}+b_4\mathbf{B^{4}})^\mu,
\end{split}
\end{equation}
with the meaning \eqref{not}.

\begin{prop}\label{d m wig}
	Fix $u\in \mathcal{S}(\mathbb{R}^{2N})$, $\alpha\in\mathbb{N}^N_0$, and let $Wig$ be the transformation defined by \eqref{cc}. Then the following properties hold:
	\begin{itemize}
		\item[(1)] $\mathbf{D}_{\mathbf{f}}^\alpha Wig[u]=Wig[(\mathbf{D_f}+\mathbf{D_s})^\alpha u]$;
		\item[(2)] $\mathbf{D}_{\mathbf{s}}^\alpha Wig[u]= Wig[(\mathbf{M_s}-\mathbf{M_f})^\alpha u]$;
		\item[(3)] $\mathbf{M}_{\mathbf{f}}^\alpha Wig[u]=Wig[\left(\frac{\mathbf{M_s}+\mathbf{M_f}}{2}\right)^\alpha u]$;
		\item[(4)] $\mathbf{M}_{\mathbf{s}}^\alpha Wig[u]=Wig\left[\left(\frac{\mathbf{D_f}-\mathbf{D_s}}{2}\right)^\alpha u\right]$.
	\end{itemize}
\end{prop}
\begin{pf}
We first observe that both $\mathbf{D_f}+\mathbf{D_s}$, $\mathbf{M_s}-\mathbf{M_f}$, $\frac{\mathbf{M_s}+\mathbf{M_f}}{2}$, and $\frac{\mathbf{D_f}-\mathbf{D_s}}{2}$ satisfy the commutation relation \eqref{commutGeneral}, so it is enough to prove the thesis for $|\alpha|=1$ and then the general case follows by \eqref{AddSumOper}.
\begin{itemize}
\item[(1)]
Let $\alpha=e_j$, where $e_j$ is the $j$-th vector of the canonical basis of $\mathbb{R}^N$. We can differentiate under the integral sign, obtaining
	\begin{align*}
		\mathbf{D}_{\mathbf{f}}^{e_j} Wig[u](x,\xi)&= D_{x_j} \int_{\mathbb{R}^N}\exp(-it\xi)u\left(x+\frac{t}{2},x-\frac{t}{2}\right)\, dt \\ &= \int_{\mathbb{R}^N}\exp(-it\xi)D_{x_j} \left(u\left(x+\frac{t}{2},x-\frac{t}{2}\right)\right)\, dt = Wig[(\mathbf{D_f}+\mathbf{D_s})^{e_j}u](x,\xi).
	\end{align*}

\item[(2)]
As in the previous case, differentiating under the integral sign we get
	\begin{align*}
		\mathbf{D}_{\mathbf{s}}^{e_j} Wig[u](x,\xi)&= D_{\xi_j} \int_{\mathbb{R}^N}\exp(-it\xi)u\left(x+\frac{t}{2},x-\frac{t}{2}\right)\, dt \\ &=\int_{\mathbb{R}^N}-t_j\exp(-it\xi)u\left(x+\frac{t}{2},x-\frac{t}{2}\right)\, dt \\ &= \int_{\mathbb{R}^N}\exp(-it\xi)\left(x_j-\frac{t_j}{2}\right)u\left(x+\frac{t}{2},x-\frac{t}{2}\right)\, dt\\& -\int_{\mathbb{R}^N}\exp(-it\xi)\left(x_j+\frac{t_j}{2}\right)u\left(x+\frac{t}{2},x-\frac{t}{2}\right)\, dt = Wig[(\mathbf{M_s}-\mathbf{M_f})^{e_j}u](x,\xi).
	\end{align*}
\item[(3)]
We have
	\begin{align*}
		\mathbf{M}_{\mathbf{f}}^{e_j} Wig[u](x,\xi)&= x_j\int_{\mathbb{R}^N}\exp(-it\xi)u\left(x+\frac{t}{2},x-\frac{t}{2}\right)\, dt \\ &=\int_{\mathbb{R}^N}\exp(-it\xi)\frac{1}{2}\left(x_j-\frac{t_j}{2}\right)u\left(x+\frac{t}{2},x-\frac{t}{2}\right)\, dt\\& +\int_{\mathbb{R}^N}\exp(-it\xi)\frac{1}{2}\left(x_j+\frac{t_j}{2}\right)u\left(x+\frac{t}{2},x-\frac{t}{2}\right)\, dt \\&= Wig\left[\left(\frac{\mathbf{M_s}+\mathbf{M_f}}{2}\right)^{e_j}u\right](x,\xi).
	\end{align*}
\item[(4)]
By integration by parts we obtain
	\begin{align*}
		\mathbf{M}_{\mathbf{s}}^{e_j} Wig[u](x,\xi)&= \xi_j \int_{\mathbb{R}^N}\exp(-it\xi)u\left(x+\frac{t}{2},x-\frac{t}{2}\right)\, dt \\&= \frac{1}{2}\int_{\mathbb{R}^N}\exp(-it\xi)\mathbf{D}_{\mathbf{f}}^{e_j} u\left(x+\frac{t}{2},x-\frac{t}{2}\right)\,dt\\&-\frac{1}{2}\int_{\mathbb{R}^N}\exp(-it\xi) \mathbf{D}_{\mathbf{s}}^{e_j} u\left(x+\frac{t}{2},x-\frac{t}{2}\right)\, dt
= Wig\left[\left(\frac{\mathbf{D_f}-\mathbf{D_s}}{2}\right)^{e_j} u\right](x,\xi).
	\end{align*}\qed
\end{itemize}
\end{pf}
\begin{oss}
	From Proposition $\ref{d m wig}$ (1) and (4) we get
	\begin{eqnarray*}
		Wig[\mathbf{D}_{\mathbf{f}}^\alpha u]=\left(\mathbf{M_s}+\frac{\mathbf{D_f}}{2}\right)^\alpha Wig[u], \\
		Wig[\mathbf{D}_{\mathbf{s}}^\alpha u]=\left(\frac{\mathbf{D_f}}{2}-\mathbf{M_s}\right)^\alpha Wig[u];
	\end{eqnarray*}
	analogously, from Proposition $\ref{d m wig}$ (2) and (3) we get
	\begin{eqnarray*}
		Wig[\mathbf{M}_{\mathbf{f}}^\alpha u]=\left(\mathbf{M_f}-\frac{\mathbf{D_s}}{2}\right)^\alpha Wig[u], \\
		Wig[\mathbf{M}_{\mathbf{s}}^\alpha u]=\left(\mathbf{M_f}+\frac{\mathbf{D_s}}{2}\right)^\alpha Wig[u],
	\end{eqnarray*}
for every $\alpha\in\mathbb{N}^N_0$. Indeed, for $\alpha=e_j$, $j=1,\dots,N$, such formulas are an easy consequence of Proposition \ref{d m wig}, and the general case follows from \eqref{AddSumOper}, since all the operators in consideration satisfy the commutation relation \eqref{commutGeneral}.
\end{oss}

\begin{oss}\label{RemCommut}
	Note that for every $\alpha,\beta\in\mathbb{N}^N_0$ and for every $u\in\mathcal{S}(\mathbb{R}^{2N})$ we have
\begin{equation}\label{commSum1}
(\mathbf{D_f}+\mathbf{D_s})^\alpha (\mathbf{M_s}-\mathbf{M_f})^\beta u=(\mathbf{M_s}-\mathbf{M_f})^\beta (\mathbf{D_f}+\mathbf{D_s})^\alpha u
\end{equation}
and
\begin{equation}\label{commSum2}
(\mathbf{D_f}-\mathbf{D_s})^\alpha (\mathbf{M_s}+\mathbf{M_f})^\beta u=(\mathbf{M_s}+\mathbf{M_f})^\beta (\mathbf{D_f}-\mathbf{D_s})^\alpha u.
\end{equation}
Indeed, from Proposition \ref{d m wig} (1), (2), we have
\begin{eqnarray*}
Wig[(\mathbf{D_f}+\mathbf{D_s})^\alpha (\mathbf{M_s}-\mathbf{M_f})^\beta u] &=& \mathbf{D}_{\mathbf{s}}^\beta \mathbf{D}_{\mathbf{f}}^\alpha Wig[u] = \mathbf{D}_{\mathbf{f}}^\alpha \mathbf{D}_{\mathbf{s}}^\beta Wig[u] \\
&=& Wig[(\mathbf{M_s}-\mathbf{M_f})^\beta (\mathbf{D_f}+\mathbf{D_s})^\alpha u],
\end{eqnarray*}
and applying $Wig^{-1}$ we get \eqref{commSum1}. Analogously, \eqref{commSum2} follows from Proposition \ref{d m wig} (3), (4).
\end{oss}

Given a linear partial differential operator with polynomial coefficients as in \eqref{ante} we denote
\begin{equation*}
\overline{P}(\mathbf{M_f},\mathbf{M_s},\mathbf{D_f},\mathbf{D_s})=P\left(\frac{\mathbf{M_s}+\mathbf{M_f}}{2}, \frac{\mathbf{D_f}-\mathbf{D_s}}{2},\mathbf{D_f}+\mathbf{D_s},\mathbf{M_s}-\mathbf{M_f}\right),
\end{equation*}
with the meaning \eqref{not1}.
\begin{prop}\label{p e p s}
	Let $P(x,y,D_x,D_y)$ be a linear partial differential operator as in \eqref{ante}. Then for each $u\in \mathcal{S}(\mathbb{R}^{2N})$, the following formula holds:
	\begin{equation}
	P(\mathbf{M_f},\mathbf{M_s},\mathbf{D_f},\mathbf{D_s}) Wig[u]=Wig\left[\overline{P}(\mathbf{M_f},\mathbf{M_s},\mathbf{D_f},\mathbf{D_s})u\right].
	\end{equation}
\end{prop}
\begin{pf}
	From Proposition $\ref{d m wig}$ we have that for each $\alpha,\beta,\gamma,\mu \in \mathbb{N}_0^N$
	\begin{eqnarray*}
		\mathbf{M}_{\mathbf{f}}^\alpha \mathbf{M}_{\mathbf{s}}^\beta \mathbf{D}_{\mathbf{f}}^\gamma \mathbf{D}_{\mathbf{s}}^\mu Wig[u]
		=Wig\left[\left(\frac{\mathbf{M_s}+\mathbf{M_f}}{2}\right)^\alpha \left(\frac{\mathbf{D_f}-\mathbf{D_s}}{2}\right)^\beta \left(\mathbf{D_f}+\mathbf{D_s}\right)^\gamma \left(\mathbf{M_s}-\mathbf{M_f}\right)^\mu u\right]
	\end{eqnarray*}
	and hence the thesis, since $Wig$ is linear. \qed
\end{pf}
Analogously if we denote
\begin{equation*}
\widetilde{P}(\mathbf{M_f},\mathbf{M_s},\mathbf{D_f},\mathbf{D_s})= P\left(\mathbf{M_f}-\frac{\mathbf{D_s}}{2},\mathbf{M_f}+\frac{\mathbf{D_s}}{2}, \mathbf{M_s}+\frac{\mathbf{D_f}}{2},\frac{\mathbf{D_f}}{2}-\mathbf{M_s}\right),
\end{equation*}
with the same scheme of the previous proof we can show the following result.
\begin{prop}\label{p e p s 2}
	Let $P(x,y,D_x,D_y)$ be a linear partial differential operator with polynomial coefficients. Then for each $u\in \mathcal{S}(\mathbb{R}^{2N})$, the following formula holds:
	\begin{equation}
	Wig\left[P\left(\mathbf{M_f},\mathbf{M_s},\mathbf{D_f},\mathbf{D_s}\right)u\right]= \widetilde{P}(\mathbf{M_f},\mathbf{M_s},\mathbf{D_f},\mathbf{D_s}) Wig[u].
	\end{equation}
\end{prop}
From now on, we consider a collection of weight functions $\omega_{i,j}$, $\sigma_{i,j}$ for $i\in\{1,2\}$ and $j=1,\dots,N$, and we fix the weights $\Omega$, $\Sigma$, $\Omega_1$ and $\Sigma_1$ in $\mathbb{R}^{2N}$ as in (\ref{Omega}), (\ref{Sigma}), (\ref{Omega1}) and (\ref{Sigma1}). Let $P(x,y,D_x,D_y)$ be a linear partial differential operator as in $(\ref{ante})$. As in the classical Schwartz case, it is easy to show that
\begin{eqnarray*}
	&&P : \mathcal{S}^{\Sigma}_{\Omega}(\mathbb{R}^{2N})\to \mathcal{S}^{\Sigma}_{\Omega}(\mathbb{R}^{2N}) \\
	&&P : (\mathcal{S}^{\Sigma}_{\Omega})'(\mathbb{R}^{2N})\to (\mathcal{S}^{\Sigma}_{\Omega})'(\mathbb{R}^{2N}).
\end{eqnarray*}
\begin{defn}
	Let $P(x,y,D_x,D_y)$ be a linear partial differential operator with polynomial coefficients and let $\Omega$, $\Sigma$ be weight functions. We say that $P$ is $\mathcal{S}^{\Sigma}_{\Omega}$-regular if
	\begin{equation}
	Pu \in \mathcal{S}^{\Sigma}_{\Omega}(\mathbb{R}^{2N}) \implies u \in \mathcal{S}^{\Sigma}_{\Omega}(\mathbb{R}^{2N}), \, \text{for each}\, u \in (\mathcal{S}^{\Sigma}_{\Omega})'(\mathbb{R}^{2N}).
	\end{equation}
\end{defn}

\begin{thm}\label{reg}
	Let $P(x,y,D_x,D_y)$ be a linear partial differential operator with polynomial coefficients, and let $\Omega$, $\Sigma$ be weight functions. If $P$ is $\mathcal{S}^{\Sigma}_{\Omega}$-regular then $\overline{P}$ is $\mathcal{S}_{\Omega_1}^{\Sigma_1}$-regular.
\end{thm}

\begin{pf}
	Suppose that $P$ is $\mathcal{S}^{\Sigma}_{\Omega}$-regular. Fix $ u \in (\mathcal{S}_{\Omega_1}^{\Sigma_1})'(\mathbb{R}^{2N})$ and suppose that $\overline{P}u \in \mathcal{S}_{\Omega_1}^{\Sigma_1}(\mathbb{R}^{2N})$. From Proposition $\ref{p e p s}$ we know that
	\begin{equation*}
	P Wig[u]=Wig\left[\overline{P}u\right],
	\end{equation*}
	hence from Proposition \ref{sim} we have $Wig\left[\overline{P}u\right] \in \mathcal{S}^{\Sigma}_{\Omega}(\mathbb{R}^{2N})$.
	Moreover, again by Proposition \ref{sim}, $Wig[u] \in(\mathcal{S}^{\Sigma}_{\Omega})'(\mathbb{R}^{2N})$. Since $P$ is $\mathcal{S}^{\Sigma}_{\Omega}$-regular, we then get that $Wig[u] \in\mathcal{S}^{\Sigma}_{\Omega}(\mathbb{R}^{2N})$. Finally applying $Wig^{-1}$ we obtain from Proposition \ref{flo} that $ u \in \mathcal{S}_{\Omega_1}^{\Sigma_1}(\mathbb{R}^{2N})$, and so $\overline{P}$ is $\mathcal{S}_{\Omega_1}^{\Sigma_1}$-regular. \qed	
\end{pf}
\begin{thm}\label{reg2}
	Let $P(x,y,D_x,D_y)$ be a linear partial differential operator with polynomial coefficients and let $\Omega$, $\Sigma$ be weight functions. If $P$ is $\mathcal{S}^{\Sigma}_{\Omega}$-regular, then $\widetilde{P}$ is $\mathcal{S}_{\Omega_1}^{\Sigma_1}$-regular.
\end{thm}
\begin{pf}
	Suppose that $P$ is $\mathcal{S}^{\Sigma}_{\Omega}$-regular. Fix $ u \in (\mathcal{S}_{\Omega_1}^{\Sigma_1})'(\mathbb{R}^{2N})$ and suppose that $\widetilde{P}u \in \mathcal{S}_{\Omega_1}^{\Sigma_1}(\mathbb{R}^{2N})$. From Proposition $\ref{p e p s 2}$ we know that for each $w \in \mathcal{S}(\mathbb{R}^{2N})$
	\begin{equation*}
	Wig\left[Pw\right]=\widetilde{P}Wig[w].
	\end{equation*}
	Consider $w=Wig^{-1}(u)\in (\mathcal{S}^{\Sigma}_{\Omega})'(\mathbb{R}^{2N})$ from Proposition \ref{flo}. Hence we obtain that
	\begin{equation*}
	Wig\left[Pw\right]=\widetilde{P}u,
	\end{equation*}
	and so $Wig\left[Pw\right] \in \mathcal{S}_{\Omega_1}^{\Sigma_1}(\mathbb{R}^{2N})$. Applying again $Wig^{-1}$ we get from Proposition \ref{flo} that $Pw \in \mathcal{S}^{\Sigma}_{\Omega}(\mathbb{R}^{2N})$. Since  $P$ is $\mathcal{S}^{\Sigma}_{\Omega}$-regular we then obtain that $w=Wig^{-1}(u)\in \mathcal{S}^{\Sigma}_{\Omega}(\mathbb{R}^{2N})$. Applying the operator $Wig$ we get from proposition \ref{sim} that $ u \in \mathcal{S}_{\Omega_1}^{\Sigma_1}(\mathbb{R}^{2N})$, and so $\widetilde{P}$ is $\mathcal{S}_{\Omega_1}^{\Sigma_1}$-regular. \qed	
\end{pf}

Now we give some examples of applications of our results in order to find classes of regular partial differential operators with polynomial coefficients.
\begin{prop}\label{ult}
Consider a multiplication operator by a polynomial, i.e. $P(x,y,D_x,D_y)=p(x,y)$, for some polynomial $p$. Then
	$P$ is $\mathcal{S}_\Omega^\Sigma$-regular if and only if $p(x,y)\neq 0$.
\end{prop}
\begin{pf}
	Suppose that there exists a point $(\overline{x},\overline{y})\in \mathbb{R}^{2N}$ such that $p(\overline{x},\overline{y})=0$. Then $p \cdot \delta_{(\overline{x},\overline{y})}=0$, i.e. $P\delta_{(\overline{x},\overline{y})}=0$. Therefore, $P$ is not $\mathcal{S}_\Omega^\Sigma$-regular. In the opposite direction, fix  $u\in(\mathcal{S}_{\Omega}^{\Sigma})'(\mathbb{R}^{2N})$ and suppose that $Pu \in \mathcal{S}_{\Omega}^{\Sigma}(\mathbb{R}^{2N})$. Hence there exists $w\in\mathcal{S}_{\Omega}^{\Sigma}(\mathbb{R}^{2N})$ such that $p(x,y)\cdot u(x,y)=w(x,y)$. Since $p(x,y)\neq 0$, it is easy to show that $u(x,y)=\frac{w(x,y)}{p(x,y)} \in \mathcal{S}_{\Omega}^{\Sigma}(\mathbb{R}^{2N})$. \qed
\end{pf}

\begin{prop}\label{ultt}
Consider an operator with constant coefficients, i.e. $P(x,y,D_x,D_y)=P(D_x,D_y)$, with symbol the polynomial $p(\xi,\eta)$. Then
	$P$ is $\mathcal{S}_\Omega^\Sigma$-regular if and only if $p(\xi,\eta)\neq 0$.
\end{prop}
\begin{pf}
	We have that $Pu=\mathcal{F}^{-1}(p\cdot\mathcal{F}u)$. Using the previous Proposition $\ref{ult}$ and Remark \ref{s} we get the claim. \qed
\end{pf}
We then obtain the following result.
\begin{corollario}\label{examples}
	Let $p(z,\zeta)=\sum_{|\alpha+\beta|\leq m} c_{\alpha\beta}z^\alpha\zeta^\beta$, for $c_{\alpha\beta}\in \mathbb{C}$ and $z,\zeta\in\mathbb{R}^N$ a polynomial in $\mathbb{R}^{2N}$, with $p(z,\zeta)\neq 0$ for every $(z,\zeta)\in \mathbb{R}^{2N}$. Then the following operators are $\mathcal{S}^{\Sigma}_{\Omega}$-regular:
	\begin{eqnarray*}
		P_1 &=& \sum_{|\alpha+\beta|\leq m} c_{\alpha\beta} \left(\frac{x+y}{2}\right)^\alpha \left(\frac{D_x-D_y}{2}\right)^\beta, \\
		P_2 &=& \sum_{|\alpha+\beta|\leq m} c_{\alpha\beta} (D_x+D_y)^\alpha (y-x)^\beta, \\
		P_3 &=& \sum_{|\alpha+\beta|\leq m} c_{\alpha\beta} \left(x-\frac{D_y}{2}\right)^\alpha \left(x+\frac{D_y}{2}\right)^\beta,\\
		P_4 &=& \sum_{|\alpha+\beta|\leq m} c_{\alpha\beta} \left(y+\frac{D_x}{2}\right)^\alpha \left(\frac{D_x}{2}-y\right)^\beta.
	\end{eqnarray*}
\end{corollario}
\begin{pf}
	It is an immediate consequence of Propositions \ref{ult}, \ref{ultt}, and Theorems \ref{reg}, \ref{reg2}. \qed
\end{pf}

\section{Regularity of time-frequency representations in the Cohen class with kernel in $\mathcal{S}'$}

In this section we extend the results of the preceding section to the case of representations in the Cohen class based on the transformation \eqref{cc}.
\begin{defn}
	Given a kernel $\kappa\in \mathcal{S}'(\mathbb{R}^{2N})$, we define the operator $Q[u]=\kappa\star Wig[u]$, for each $u\in \mathcal{S}(\mathbb{R}^{2N})$. An operator $Q$ is called a time-frequency representation in the Cohen class with kernel $\kappa$.
\end{defn}

We shall consider in particular kernels $\kappa(x,y)$ whose Fourier transform $\hat{\kappa}(\xi,\eta)$ is given by
\begin{equation}\label{kernel}
\hat{\kappa}(\xi,\eta)=\exp\left(-i\sum_{j=1}^N p_j(\xi_j,\eta_j)\right),
\end{equation}
where $p_j$, $j=1,\dots,N$ are polynomials in $\mathbb{R}^2$ of any order, with coefficients in $\mathbb{R}$.

Now consider the polynomial $p_j=p_j(s,t)$, for $s,t\in\mathbb{R}$; we indicate with
$$
(\partial_1 p_j)(s,t) = \frac{\partial p_j}{\partial s}(s,t),\quad (\partial_2 p_j)(s,t) = \frac{\partial p_j}{\partial t}(s,t).
$$
Then, for $j=1,\dots,N$ we define the operators
\begin{eqnarray}
&&R_j=R_j(D_x,D_y)=(\partial_1 p_j)(D_{x_j},D_{y_j}), \label{Rj}\\
&&T_j=T_j(D_x,D_y)=(\partial_2 p_j)(D_{x_j},D_{y_j}), \label{Tj}\\
&&R_j^*=R_j^*(D_x,D_y)=(\partial_1 p_j)(D_{x_j}+D_{y_j},y_j-x_j), \notag\\
&&T_j^*=T_j^*(D_x,D_y)=(\partial_2 p_j)(D_{x_j}+D_{y_j},y_j-x_j), \notag
\end{eqnarray}
for $x,y\in\mathbb{R}^N$. Observe that $R_j$ and $T_j$, are partial differential operators with constant coefficients in $\mathbb{R}^{2N}$, while $R_j^*$ and $T_j^*$ are partial differential operators with polynomial coefficients in $\mathbb{R}^{2N}$. We then define the corresponding vectors of operators
\begin{eqnarray}
&&\mathbf{R}=(R_1,\dots,R_N),\quad \mathbf{T}=(T_1,\dots,T_N), \label{RT}\\
&&\mathbf{R^*}=(R_1^*,\dots,R_N^*),\quad \mathbf{T^*}=(T_1^*,\dots,T_N^*). \label{RstarTstar}
\end{eqnarray}

\begin{oss}\label{RWig}
By Proposition \ref{p e p s} we have that for every $u\in\mathcal{S}(\mathbb{R}^N)$
$$
\mathbf{R}Wig[u] = Wig[\mathbf{R^*} u],\ \ \text{and}\ \ \mathbf{T}Wig[u] = Wig[\mathbf{T^*} u],
$$
where $Wig$ applied to a vector is intended as the vector of $Wig$ applied to the components, so that the equalities above are intended componentwise. Then, since $\mathbf{R}$ and $\mathbf{T}$ are vectors of operators with constant coefficients, by the properties of the convolution we also have
$$
\mathbf{R}Q[u] = \kappa * \mathbf{R}Wig[u] = \kappa * Wig[\mathbf{R^*} u] = Q[\mathbf{R^*} u],
$$
and analogously
$$
\mathbf{T} Q[u] = Q[\mathbf{T^*} u].
$$

\end{oss}

\begin{lemma}\label{lemmAgg}
Let $\mathbf{R}$, $\mathbf{T}$ as before, $\mathbf{M_f}$, $\mathbf{M_s}$ as in Definition \ref{defMD}, and $\kappa$ as in \eqref{kernel}. Then
$$
\mathbf{M_f} \kappa=\mathbf{R} \kappa\ \ \text{and}\ \ \mathbf{M_s} \kappa=\mathbf{T} \kappa.
$$
\end{lemma}

\begin{pf}
In order to prove the equality $\mathbf{M_f} \kappa=\mathbf{R} \kappa$ we have to show that for every $j=1,\dots,N$,
$$
x_j\kappa(x,y)=(\partial_1 p_j)(D_{x_j},D_{y_j}) \kappa (x,y).
$$
Applying the Fourier transform to both sides, this in turn is equivalent to
$$
i\partial_{\xi_j} \hat{\kappa}(\xi,\eta) = (\partial_1 p_j)(\xi_j,\eta_j) \hat{\kappa}(\xi,\eta),
$$
and this last equality trivially follows from \eqref{kernel}. The other relation $\mathbf{M_s} \kappa=\mathbf{T} \kappa$ can be proved in the same way. \qed
\end{pf}

\begin{prop}\label{d m q}
	Fix $u\in \mathcal{S}(\mathbb{R}^{2N})$, and let $\kappa$ be the kernel defined by \eqref{kernel}. Then for every $\alpha\in\mathbb{N}^N_0$ the following properties hold:
	\begin{itemize}
		\item[(1)] $\mathbf{D}_{\mathbf{f}}^\alpha Q[u]=Q[(\mathbf{D_f}+\mathbf{D_s})^\alpha u]$;
		\item[(2)] $\mathbf{D}_{\mathbf{s}}^\alpha Q[u]= Q[(\mathbf{M_s}-\mathbf{M_f})^\alpha u]$;
		\item[(3)] $\mathbf{M}_{\mathbf{f}}^\alpha Q[u]=Q[\left(\frac{\mathbf{M_s}+\mathbf{M_f}}{2}+\mathbf{R^*}\right)^\alpha u]$;
		\item[(4)] $\mathbf{M}_{\mathbf{s}}^\alpha Q[u]=Q\left[\left(\frac{\mathbf{D_f}-\mathbf{D_s}}{2}+\mathbf{T^*}\right)^\alpha u\right]$.
	\end{itemize}
\end{prop}
\begin{pf}
	In the following the integrals are intended as the action of the distribution $\kappa$ when $\kappa$ is not a function. \\[0.1cm]
(1) and (2) are trivial consequences of Proposition \ref{d m wig}, since $\mathbf{D}_{\mathbf{f}}^\alpha Q[u] = \kappa\star \mathbf{D}_{\mathbf{f}}^\alpha Wig[u]$ and $\mathbf{D}_{\mathbf{s}}^\alpha Q[u]=\kappa\star \mathbf{D}_{\mathbf{s}}^\alpha Wig[u]$. \\[0.1cm]
	Concerning point (3) we first observe that the vector of operators $\frac{\mathbf{M_s}+\mathbf{M_f}}{2}+\mathbf{R^*}$ satisfies the commutation relations \eqref{commutGeneral}, so it satisfies \eqref{AddSumOper}; this means that it is enough to prove the thesis for $|\alpha|=1$, and the general case follows by induction.

For $\alpha=e_j$, $j=1,\dots,N$, we get that
	\begin{align*}
		\mathbf{M}_{\mathbf{f}}^{e_j} Q[u](x,\xi)&= \int_{\mathbb{R}^{2N}} x_j\kappa(t,s)Wig(x-t,\xi-s)\, dtds \\&
=\int_{\mathbb{R}^{2N}} (x_j-t_j)\kappa(t,s)Wig(x-t,\xi-s)\, dtds \\&+\int_{\mathbb{R}^{2N}} t_j\kappa(t,s)Wig(x-t,\xi-s)\, dtds \\
&= \kappa\star (\mathbf{M}_{\mathbf{f}}^{e_j}Wig[u])+(\mathbf{M}_{\mathbf{f}}^{e_j} \kappa)\star Wig[u].
	\end{align*}
Now we observe that, by Lemma \ref{lemmAgg}, $\mathbf{M}_{\mathbf{f}}^{e_j} \kappa = \mathbf{R}^{e_j} \kappa$, and since $\mathbf{R}^{e_j}$ is an operator with constant coefficients we have that $(\mathbf{M}_{\mathbf{f}}^{e_j} \kappa)\star Wig[u] = \kappa\star \mathbf{R}^{e_j} Wig[u]$. Then, applying Proposition \ref{d m wig} and Remark \ref{RWig} we get
	\begin{align*}
		\mathbf{M}_{\mathbf{f}}^{e_j} Q[u](x,\xi)
		&=\kappa\star Wig\left[\left(\frac{\mathbf{M_s}+\mathbf{M_f}}{2}\right)^{e_j} u\right]+\kappa\star Wig[(\mathbf{R^*})^{e_j}u] \\&
=Q\left[\left(\frac{\mathbf{M_s}+\mathbf{M_f}}{2}+\mathbf{R^*}\right)^{e_j} u\right].
	\end{align*}
    To prove (4) we first observe that the vector of operators $\frac{\mathbf{D_f}-\mathbf{D_s}}{2}+\mathbf{T^*}$ satisfies the commutation relations \eqref{commutGeneral}, so it satisfies \eqref{AddSumOper}; as before it is enough to prove the thesis for $|\alpha|=1$, and the general case follows by induction. For $\alpha=e_j$, $j=1,\dots,N$, we get that
    \begin{align*}
    	\mathbf{M}_{\mathbf{s}}^{e_j} Q[u](x,\xi)&= \int_{\mathbb{R}^{2N}} \xi_j \kappa(t,s)Wig(x-t,\xi-s)\, dtds \\&
    =\int_{\mathbb{R}^{2N}} (\xi_j-s_j)\kappa(t,s)Wig(x-t,\xi-s)\, dtds \\&+\int_{\mathbb{R}^{2N}} s_j \kappa(t,s)Wig(x-t,\xi-s)\, dtds \\&
    = \kappa\star (\mathbf{M}_{\mathbf{s}}^{e_j} Wig[u])+(\mathbf{M}_{\mathbf{s}}^{e_j} \kappa)\star Wig[u].
    \end{align*}
By Lemma \ref{lemmAgg} we have $\mathbf{M}_{\mathbf{s}}^{e_j} \kappa = \mathbf{T}^{e_j} \kappa$, and since $\mathbf{T}^{e_j}$ is an operator with constant coefficients we have that $(\mathbf{M}_{\mathbf{s}}^{e_j} \kappa)\star Wig[u] = \kappa\star \mathbf{T}^{e_j} Wig[u]$. Then, applying Proposition \ref{d m wig} and Remark \ref{RWig} we get
    \begin{align*}
    	\mathbf{M}_{\mathbf{s}}^{e_j} Q[u](x,\xi)
    	&=\kappa\star Wig\left[\left(\frac{\mathbf{D_f}-\mathbf{D_s}}{2}\right)^{e_j}u\right]+\kappa\star Wig[(\mathbf{T^*})^{e_j} u] \\&
    = Q\left[\left(\frac{\mathbf{D_f}-\mathbf{D_s}}{2}+\mathbf{T^*}\right)^{e_j}u\right].
    \end{align*}
    \qed
\end{pf}

\begin{oss}\label{RemQDM}
The following relations hold:
\begin{eqnarray}
&&Q[\mathbf{M}_{\mathbf{f}}^\alpha u] = \left(\mathbf{M_f}-\frac{\mathbf{D_s}}{2}-\mathbf{R}\right)^\alpha Q[u], \quad
Q[\mathbf{M}_{\mathbf{s}}^\alpha u] = \left(\mathbf{M_f}+\frac{\mathbf{D_s}}{2}-\mathbf{R}\right)^\alpha Q[u], \label{Qdermult1} \\
&&Q[\mathbf{D}_{\mathbf{f}}^\alpha u] = \left(\frac{\mathbf{D}_{\mathbf{f}}}{2} + \mathbf{M}_{\mathbf{s}} - \mathbf{T}\right)^\alpha Q[u], \quad
Q[\mathbf{D}_{\mathbf{s}}^\alpha u] = \left(\frac{\mathbf{D_f}}{2}-\mathbf{M_s}+\mathbf{T}\right)^\alpha Q[u]; \label{Qdermult2}
\end{eqnarray}
indeed, from Proposition $\ref{d m q}$ with $\alpha=e_j$, $j=1,\dots,N$, Remark \ref{RWig}, and the linearity of $Q$, we get
	\begin{eqnarray*}
    &&\mathbf{D}_{\mathbf{f}}^{e_j} Q[u]=Q[\mathbf{D}_{\mathbf{f}}^{e_j}u] + Q[\mathbf{D}_{\mathbf{s}}^{e_j}u], \\
    &&\mathbf{D}_{\mathbf{s}}^{e_j} Q[u]=Q[\mathbf{M}_{\mathbf{s}}^{e_j}u] - Q[\mathbf{M}_{\mathbf{f}}^{e_j}u], \\
    &&\mathbf{M}_{\mathbf{f}}^{e_j} Q[u]=\frac{1}{2}Q[\mathbf{M}_{\mathbf{s}}^{e_j}u] + \frac{1}{2}Q[\mathbf{M}_{\mathbf{f}}^{e_j}u] + \mathbf{R}^{e_j} Q[u], \\
    &&\mathbf{M}_{\mathbf{s}}^{e_j} Q[u]=\frac{1}{2}Q[\mathbf{D}_{\mathbf{f}}^{e_j}u] - \frac{1}{2}Q[\mathbf{D}_{\mathbf{s}}^{e_j}u] + \mathbf{T}^{e_j} Q[u].
	\end{eqnarray*}
Then, combining these last relations we easily get \eqref{Qdermult1} and \eqref{Qdermult2} for $\alpha=e_j$; the general case follows from the fact that the four vectors of operators $(\mathbf{M_f}-\frac{\mathbf{D_s}}{2}-\mathbf{R})$, $(\mathbf{M_f}+\frac{\mathbf{D_s}}{2}-\mathbf{R})$, $(\frac{\mathbf{D}_{\mathbf{f}}}{2} - \mathbf{M}_{\mathbf{s}} - \mathbf{T})$ and $(\frac{\mathbf{D_f}}{2}-\mathbf{M_s}+\mathbf{T})$ satisfy the commutation relations \eqref{commutGeneral}, so they satisfy \eqref{AddSumOper}.
\end{oss}

\begin{oss}\label{notanucleo}
The reason for the particular choice of $\kappa$ as in \eqref{kernel}, with each polynomial depending only on a single couple of variables $(\xi_j,\eta_j)$ is related to Proposition \ref{d m q} and Remark \ref{RemQDM}, since this choice ensures us that all the vectors of operators that we consider in those results satisfy the commutation relation \eqref{commutGeneral}. It would be interesting to analyze the case when $\kappa$ is chosen in such a way that $\hat{\kappa}(\xi,\eta)=\exp(-ip(\xi,\eta))$ for a generic polynomial $p$ in $\mathbb{R}^{2N}$ with real coefficients. However, in this case it looks very difficult to have good formulas corresponding to the ones of Proposition \ref{d m q}.
\end{oss}

\begin{prop}\label{b q}
Let $B(x,y,D_x,D_y)$ be a linear partial differential operator with polynomial coefficients and let the kernel $\kappa\in \mathcal{S}'(\mathbb{R}^{2N})$ be defined by \eqref{kernel}. Then for each $u\in \mathcal{S}(\mathbb{R}^{2N})$, the time-frequency representation $Q[w]=\kappa\star Wig[w]$ satisfies:
\begin{equation}
B(\mathbf{M_f},\mathbf{M_s},\mathbf{D_f},\mathbf{D_s}) Q[u]=Q\left[\overline{B}(\mathbf{M_f},\mathbf{M_s}, \mathbf{D_f},\mathbf{D_s})u\right],
\end{equation}
where $\overline{B}$ is the linear partial differential operator with polynomial coefficients defined by
\begin{equation}
\overline{B}(\mathbf{M_f},\mathbf{M_s},\mathbf{D_f},\mathbf{D_f}) =B\left(\frac{\mathbf{M_s}+\mathbf{M_f}}{2}+\mathbf{R^*},\frac{\mathbf{D_f}-\mathbf{D_s}}{2}+\mathbf{T^*}, \mathbf{D_f}+\mathbf{D_s}, \mathbf{M_s}-\mathbf{M_f}\right).
\end{equation}
\end{prop}

\begin{pf}
The proof follows immediately from Proposition \ref{d m q} and the linearity of $Q$. \qed
\end{pf}

\begin{prop}\label{b q 2}
	Let $B(x,y,D_x,D_y)$ be a linear partial differential operator with polynomial coefficients and let the kernel $\kappa\in \mathcal{S}'(\mathbb{R}^{2N})$ be defined by \eqref{kernel}. Then for each $u\in \mathcal{S}(\mathbb{R}^{2N})$, the time-frequency representation $Q[w]=\kappa\star Wig[w]$ satisfies:
	\begin{equation}
	Q[B(\mathbf{M_f},\mathbf{M_s},\mathbf{D_f},\mathbf{D_s})u]= \widetilde{B}(\mathbf{M_f},\mathbf{M_s},\mathbf{D_f},\mathbf{D_s})Q[u],
	\end{equation}
	where $\widetilde{B}$ is the linear partial differential operator with polynomial coefficients defined by
	\begin{equation}\label{simm}
	\widetilde{B}(\mathbf{M_f},\mathbf{M_s},\mathbf{D_f},\mathbf{D_s}) =B\left(\mathbf{M_f}-\frac{\mathbf{D_s}}{2}-\mathbf{R}, \mathbf{M_f}+\frac{\mathbf{D_s}}{2}-\mathbf{R},\frac{\mathbf{D_f}}{2}+\mathbf{M_s}-\mathbf{T}, \frac{\mathbf{D_f}}{2}-\mathbf{M_s}+\mathbf{T}\right),
	\end{equation}
\end{prop}

\begin{pf}
The proof follows immediately from Remark \ref{RemQDM} and the linearity of $Q$. \qed
\end{pf}

Now we want to study the action of $Q$ on the spaces $\mathcal{S}^\Sigma_\Omega(\mathbb{R}^{2N})$.
\begin{thm}\label{q}
	Fix the kernel $\kappa\in \mathcal{S}'(\mathbb{R}^{2N})$ as in \eqref{kernel}. Then the following properties hold:
	\begin{itemize}
		\item[(1)] $Q: \mathcal{S}'(\mathbb{R}^{2N})\to \mathcal{S}'(\mathbb{R}^{2N})$,
		\item[(2)] $Q: \mathcal{S}(\mathbb{R}^{2N})\to \mathcal{S}(\mathbb{R}^{2N})$,
		\item[(3)] $Q: (\mathcal{S}^{\Sigma}_{\Omega})'(\mathbb{R}^{2N})\to (\mathcal{S}_{\Omega_1}^{\Sigma_1})'(\mathbb{R}^{2N})$,
 		\item[(4)] $Q: \mathcal{S}^{\Sigma}_{\Omega}(\mathbb{R}^{2N})\to \mathcal{S}_{\Omega_1}^{\Sigma_1}(\mathbb{R}^{2N})$,
	\end{itemize}
and in all the cases $Q$ is invertible.
\end{thm}
\begin{pf}
	The proof of (1) and (2) can be found in \cite{8}, in the case $N=1$, and the proof works in the same way in higher dimension.
	\\[0.2cm]
	(3) Fix $u\in (\mathcal{S}^{\Sigma}_{\Omega})'(\mathbb{R}^{2N})$. We observe that $\mathcal{F}(Q[u])=\hat{\kappa}\mathcal{F}(Wig[u]) \in (\mathcal{S}^{\Omega_1}_{\Sigma_1})'(\mathbb{R}^{2N})$, since $\hat{\kappa}$ has a polynomial growth and $\mathcal{F}(Wig[u]) \in (\mathcal{S}^{\Omega_1}_{\Sigma_1})'(\mathbb{R}^{2N})$ from Proposition \ref{sim} and Remark \ref{s} extended to ultradistributions. Applying $\mathcal{F}^{-1}$ we get that $Q[u]\in (\mathcal{S}_{\Omega_1}^{\Sigma_1})'(\mathbb{R}^{2N})$ as desired.
	\\[0.2cm]
	Now we prove the invertibility of $Q$. The injectivity follows from the injectivity of the Wigner-like transform. To prove the surjectivity, fix $w\in(\mathcal{S}_{\Omega_1}^{\Sigma_1})'(\mathbb{R}^{2N})$. Then $\hat{w}\in (\mathcal{S}^{\Omega_1}_{\Sigma_1})'(\mathbb{R}^{2N})$. Since $1/\hat{\kappa}$ has still a polynomial growth, then also $\hat{w}/\hat{\kappa}\in (\mathcal{S}^{\Omega_1}_{\Sigma_1})'(\mathbb{R}^{2N})$. By the surjectivity of the Fourier transform there exists $v\in (\mathcal{S}_{\Omega_1}^{\Sigma_1})'(\mathbb{R}^{2N})$ such that $\hat{v}=\hat{w}/\hat{\kappa}$. By the surjectivity of the Wigner-like transform, $v=Wig[u]$ for some $u\in (\mathcal{S}^{\Sigma}_{\Omega})'(\mathbb{R}^{2N})$ and therefore
	\begin{equation*}
		\hat{w}=\hat{\kappa}\hat{v}=\hat{\kappa} \mathcal{F}(Wig[u])= \mathcal{F}(\kappa\star Wig[u])=\mathcal{F}(Q[u])
	\end{equation*}
	and by the injectivity of the Fourier transform $w=Q[u]$, for $u\in (\mathcal{S}^{\Sigma}_{\Omega})'(\mathbb{R}^{2N})$.
	\\[0.2cm]
	(4) Fix $u\in \mathcal{S}^{\Sigma}_{\Omega}(\mathbb{R}^{2N})$. We have that
	\begin{equation*}
		Q[u]=\kappa\star Wig[u]=\mathcal{F}^{-1}(\hat{\kappa}\mathcal{F}(Wig[u]))\in \mathcal{S}_{\Omega_1}^{\Sigma_1}(\mathbb{R}^{2N})
	\end{equation*}
	since $\hat{\kappa}$ has a polynomial growth and $\mathcal{F}(Wig[u]) \in \mathcal{S}^{\Omega_1}_{\Sigma_1}(\mathbb{R}^{2N})$ for $u\in \mathcal{S}^{\Sigma}_{\Omega}(\mathbb{R}^{2N})$. The invertibility can be proved as in the previous point. \qed
\end{pf}
\begin{oss}\label{reminvert}
By the invertibility properties of Theorem \ref{q} we have that if $u\in(\mathcal{S}_\Omega^\Sigma)'(\mathbb{R}^{2N})$ and $Q[u]\in\mathcal{S}_{\Omega_1}^{\Sigma_1}(\mathbb{R}^{2N})$, then $u\in\mathcal{S}_\Omega^\Sigma(\mathbb{R}^{2N})$.
\end{oss}

\begin{thm}
	Let $B(x,y,D_x,D_y)$ be a linear partial differential operator with polynomial coefficients, and fix the kernel $\kappa\in \mathcal{S}'(\mathbb{R}^{2N})$ as in \eqref{kernel}. If $B$ is $\mathcal{S}^{\Sigma}_{\Omega}$-regular then $\overline{B}$ is $\mathcal{S}_{\Omega_1}^{\Sigma_1}$-regular.
\end{thm}
\begin{pf}
	Suppose that $B$ is $\mathcal{S}^{\Sigma}_{\Omega}$-regular. Fix $u \in (\mathcal{S}_{\Omega_1}^{\Sigma_1})'(\mathbb{R}^{2N})$ and suppose that $\overline{B}u \in \mathcal{S}_{\Omega_1}^{\Sigma_1}(\mathbb{R}^{2N})$. From Proposition $\ref{b q}$ we know that
	\begin{equation*}
	B Q[u]=Q\left[\overline{B}u\right],
	\end{equation*}
	hence from Theorem \ref{q} we get $Q\left[\overline{B}u\right] \in \mathcal{S}^{\Sigma}_{\Omega}(\mathbb{R}^{2N})$.
	Moreover, again by Theorem \ref{q}, we have $Q[u] \in(\mathcal{S}^{\Sigma}_{\Omega})'(\mathbb{R}^{2N})$. Since $B$ is $\mathcal{S}^{\Sigma}_{\Omega}$-regular, we then get that $Q[u] \in\mathcal{S}^{\Sigma}_{\Omega}(\mathbb{R}^{2N})$. Finally by Remark \ref{reminvert} we obtain that $u \in \mathcal{S}_{\Omega_1}^{\Sigma_1}(\mathbb{R}^{2N})$, and so $\overline{B}$ is $\mathcal{S}_{\Omega_1}^{\Sigma_1}$-regular. \qed	
\end{pf}
\begin{thm}
	Let $B(x,y,D_x,D_y)$ be a linear partial differential operator with polynomial coefficients, and fix the kernel $\kappa\in \mathcal{S}'(\mathbb{R}^{2N})$ as in \eqref{kernel}. If $B$ is $\mathcal{S}^{\Sigma}_{\Omega}$-regular then $\widetilde{B}$ is $\mathcal{S}_{\Omega_1}^{\Sigma_1}$-regular.
\end{thm}
\begin{pf}
	Suppose that $B$ is $\mathcal{S}^{\Sigma}_{\Omega}$-regular. Fix $u \in (\mathcal{S}_{\Omega_1}^{\Sigma_1})'(\mathbb{R}^{2N})$ and suppose that $\widetilde{B}u \in \mathcal{S}_{\Omega_1}^{\Sigma_1}(\mathbb{R}^{2N})$. From Proposition $\ref{b q 2}$ we know that for each $w \in \mathcal{S}(\mathbb{R}^{2N})$
	\begin{equation*}
	Q\left[Bw\right]=\widetilde{B}Q[w].
	\end{equation*}
	Consider $w=Q^{-1}[u]$; from Theorem \ref{q} we have $w\in (\mathcal{S}^{\Sigma}_{\Omega})'(\mathbb{R}^{2N})$. Hence we obtain that
	\begin{equation*}
	Q\left[Bw\right]=\widetilde{B}u,
	\end{equation*}
	and so $Q\left[Bw\right] \in \mathcal{S}_{\Omega_1}^{\Sigma_1}(\mathbb{R}^{2N})$. By Remark \ref{reminvert} we then have that $Bw \in \mathcal{S}^{\Sigma}_{\Omega}(\mathbb{R}^{2N})$. Since  $B$ is $\mathcal{S}^{\Sigma}_{\Omega}$-regular we then obtain that $w=Q^{-1}[u]\in \mathcal{S}^{\Sigma}_{\Omega}(\mathbb{R}^{2N})$. Applying the operator $Q$ we get from Theorem \ref{q} that $u \in \mathcal{S}_{\Omega_1}^{\Sigma_1}(\mathbb{R}^{2N})$, and so $\widetilde{B}$ is $\mathcal{S}_{\Omega_1}^{\Sigma_1}$-regular. \qed	 \end{pf}

We can give a further generalization of the last results, by taking a kernel $\kappa_1$ of the following form; let $\kappa$ be defined as in \eqref{kernel}, and let $q\in\mathbb{C}[\xi,\eta] $ be a polynomial that never vanishes on $\mathbb{R}^{2N}$. We define $\kappa_1$ by
\begin{equation}\label{kernelK1}
	\hat{\kappa}_1(\xi,\eta)=q(\xi,\eta)\hat{\kappa}(\xi,\eta).
\end{equation}
Then $\kappa_1(x,y)=q(D_x,D_y)\kappa(x,y)$ and, by Proposition \ref{p e p s}, we have
\begin{equation}\label{q1q}
	Q_1[u]=\kappa_1\star Wig[u]=\kappa\star (q(\mathbf{D_f},\mathbf{D_s})Wig[u])	=\kappa\star Wig[Au]=Q[Au],
\end{equation}
where
\begin{equation}\label{operatorA}
A(\mathbf{M_f},\mathbf{M_s},\mathbf{D_f},\mathbf{D_s})=q(\mathbf{D_f}+\mathbf{D_s},\mathbf{M_s}-\mathbf{M_f}).
\end{equation}
We can give the following result.
\begin{thm}\label{q1}
Fix the kernel $\kappa_1\in \mathcal{S}'(\mathbb{R}^{2N})$ as in \eqref{kernelK1}, where $\kappa$ is defined by \eqref{kernel}. Writing $Q_1[u]=\kappa_1\star Wig[u]$, we have the following properties:
	\begin{itemize}
		\item[(1)] $Q_1: \mathcal{S}'(\mathbb{R}^{2N})\to \mathcal{S}'(\mathbb{R}^{2N})$,
		\item[(2)] $Q_1: \mathcal{S}(\mathbb{R}^{2N})\to \mathcal{S}(\mathbb{R}^{2N})$,
		\item[(3)] $Q_1: (\mathcal{S}^{\Sigma}_{\Omega})'(\mathbb{R}^{2N})\to (\mathcal{S}_{\Omega_1}^{\Sigma_1})'(\mathbb{R}^{2N})$,
		\item[(4)] $Q_1: \mathcal{S}^{\Sigma}_{\Omega}(\mathbb{R}^{2N})\to \mathcal{S}_{\Omega_1}^{\Sigma_1}(\mathbb{R}^{2N})$,
	\end{itemize}
and in all cases $Q_1$ is invertible. Moreover, if $u\in(\mathcal{S}_\Omega^\Sigma)'(\mathbb{R}^{2N})$ and $Q_1[u]\in\mathcal{S}_{\Omega_1}^{\Sigma_1}(\mathbb{R}^{2N})$, then $u\in\mathcal{S}_\Omega^\Sigma(\mathbb{R}^{2N})$.
\end{thm}
\begin{pf}
	The proof is analogous to that of Theorem \ref{q} and Remark \ref{reminvert}, since $\hat{\kappa}_1(\xi,\eta)=q(\xi,\eta)\widehat{\kappa}(\xi,\eta)$ and $q(\xi,\eta)$ never vanishes. \qed
\end{pf}
\begin{thm}\label{ab}
Fix the kernel $\kappa_1\in \mathcal{S}'(\mathbb{R}^{2N})$ as in \eqref{kernelK1}, where $\kappa$ is defined by \eqref{kernel}. Writing $Q_1[u]=\kappa_1\star Wig[u]$ we have that the following formula holds for $u\in\mathcal{S}(\mathbb{R}^{2N})$:
	\begin{equation}\label{for}
		Q_1[Bu]=\widetilde{AB}Q[u],
	\end{equation}
	where $A$ is the operator \eqref{operatorA}, and $\widetilde{AB}$ is obtained by $AB$ as in $\eqref{simm}$. Moreover $B$ is $\mathcal{S}^{\Sigma}_{\Omega}$-regular if and only if $\widetilde{AB}$ is $\mathcal{S}_{\Omega_1}^{\Sigma_1}$-regular.
\end{thm}
\begin{pf}
	The formula $(\ref{for})$ follows from $(\ref{q1q})$ and Proposition \ref{b q 2}.
	\\[0.2cm]
	Now we suppose that $B$ is $\mathcal{S}^{\Sigma}_{\Omega}$-regular and we prove that $\widetilde{AB}$ is $\mathcal{S}_{\Omega_1}^{\Sigma_1}$-regular. Fix $u \in (\mathcal{S}_{\Omega_1}^{\Sigma_1})'(\mathbb{R}^{2N})$ and suppose that $\widetilde{AB}u \in \mathcal{S}_{\Omega_1}^{\Sigma_1}(\mathbb{R}^{2N})$. By Theorem \ref{q} there exists $w\in(\mathcal{S}^{\Sigma}_{\Omega})'(\mathbb{R}^{2N})$ such that $Q[w]=u$. By $\eqref{for}$ we have that $Q_1[Bw]=\widetilde{AB}Q[w]=\widetilde{AB}u$ and hence $Bw \in \mathcal{S}^{\Sigma}_{\Omega}(\mathbb{R}^{2N})$ from Theorem \ref{q1}. Since $B$ is $\mathcal{S}^{\Sigma}_{\Omega}$-regular, then $w\in \mathcal{S}^{\Sigma}_{\Omega}(\mathbb{R}^{2N})$. Applying $Q$ we get that $u=Q[w]\in \mathcal{S}_{\Omega_1}^{\Sigma_1}(\mathbb{R}^{2N})$, and so $\widetilde{AB}$ is $\mathcal{S}_{\Omega_1}^{\Sigma_1}$-regular.

	Reciprocally, we assume that  $\widetilde{AB}$ is $\mathcal{S}_{\Omega_1}^{\Sigma_1}$-regular and we prove that $B$ is $\mathcal{S}^{\Sigma}_{\Omega}$-regular. Fix $u \in (\mathcal{S}^{\Sigma}_{\Omega})'(\mathbb{R}^{2N})$ and suppose that $Bu \in \mathcal{S}^{\Sigma}_{\Omega}(\mathbb{R}^{2N})$. Then $Q_1[Bu]\in \mathcal{S}_{\Omega_1}^{\Sigma_1}(\mathbb{R}^{2N})$ by Theorem \ref{q1}. Using formula $\eqref{for}$ we obtain that $\widetilde{AB}Q[u]=Q_1[Bu]\in\mathcal{S}_{\Omega_1}^{\Sigma_1}(\mathbb{R}^{2N})$. Since $Q[u]\in(\mathcal{S}_{\Omega_1}^{\Sigma_1})'(\mathbb{R}^{2N})$ and $\widetilde{AB}$ is $\mathcal{S}_{\Omega_1}^{\Sigma_1}$-regular we get that $Q[u]\in\mathcal{S}_{\Omega_1}^{\Sigma_1}(\mathbb{R}^{2N})$. By Theorem \ref{q1} we conclude that $u\in\mathcal{S}^{\Sigma}_{\Omega}(\mathbb{R}^{2N})$, and so $B$ is $\mathcal{S}^{\Sigma}_{\Omega}$-regular. \qed
\end{pf}
\begin{oss}
	In the particular case $q\equiv1$ we have that $A$ is the identity, and so Theorem \ref{ab} implies that a linear partial differential operator $B$ with polynomial coefficients  is $\mathcal{S}^{\Sigma}_{\Omega}$-regular if and only if $\widetilde{B}$ is $\mathcal{S}_{\Omega_1}^{\Sigma_1}$-regular.
\end{oss}

We conclude with some examples of application of Propositions \ref{b q} and \ref{b q 2}. We have already observed in Propositions \ref{ult} and \ref{ultt} that the operator with polynomial coefficients
$$
P(x,y,D_x,D_y)=p(x,y),
$$
resp.
$$
Q(x,y,D_x,D_y)=q(D_x,D_y),
$$
is $\mathcal{S}_\Omega^\Sigma$-regular if and only if $p(x,y)$, resp. $q(\xi,\eta)$, never vanishes. If we consider, as particular case, a kernel of the form \eqref{kernel} where we assume that the polynomials $p_j$ are of the form
$$
p_j(\xi_j,\eta_j)=p_{j,1}(\xi_j)+p_{j,2}(\eta_j),
$$
for every $j=1,\dots,N$, then the operators $R_j=R_j(D_{x_j})$ and $T_j=T_j(D_{y_j})$ defined in \eqref{Rj} and \eqref{Tj} can be chosen as arbitrary differential operators with constant real coefficients, of any order and without any other assumption on their symbols. Moreover, we have
$$
R_j^*=R_j(D_{x_j}+D_{y_j})\quad\text{and}\quad T_j^*=T_j(y_j-x_j),
$$
again for arbitrary $R_j$ and $T_j$ with real coefficients. Now let
$$
p(z,\zeta)=\sum_{|\alpha+\beta|\leq m} c_{\alpha\beta}z^\alpha\zeta^\beta,\quad c_{\alpha\beta}\in\mathbb{C},\ z,\zeta\in\mathbb{R}^N,
$$
be a polynomial that never vanishes. Then, as in Corollary \ref{examples}, from Propositions \ref{b q} and \ref{b q 2}, we obtain that the following operators are $\mathcal{S}_\Omega^\Sigma$-regular:
\begin{eqnarray*}
&&P_1=\sum_{|\alpha+\beta|\leq m} c_{\alpha\beta}\left(\frac{x_1+y_1}{2}+R_1(D_{x_1}+D_{y_1})\right)^{\alpha_1}\dots \left(\frac{x_N+y_N}{2}+R_N(D_{x_N}+D_{y_N})\right)^{\alpha_N} \\
&&\qquad\qquad\qquad\left(\frac{D_{x_1}-D_{y_1}}{2}+T_1(y_1-x_1)\right)^{\beta_1} \dots\left(\frac{D_{x_N}-D_{y_N}}{2}+T_N(y_N-x_N)\right)^{\beta_N}; \\
&&
P_2=\sum_{|\alpha+\beta|\leq m} c_{\alpha\beta}\left(x_1-\frac{D_{y_1}}{2}-R_1(D_{x_1})\right)^{\alpha_1}\dots \left(x_N-\frac{D_{y_N}}{2}-R_N(D_{x_N})\right)^{\alpha_N} \\
&&\qquad\qquad\qquad\left(x_1+\frac{D_{y_1}}{2}-R_1(D_{x_1})\right)^{\beta_1}\dots \left(x_N+\frac{D_{y_N}}{2}-R_N(D_{x_N})\right)^{\beta_N}; \\
&&P_3=\sum_{|\alpha+\beta|\leq m} c_{\alpha\beta}\left(\frac{D_{x_1}}{2}+y_1-T_1(D_{y_1})\right)^{\alpha_1}\dots \left(\frac{D_{x_N}}{2}+y_N-T_N(D_{y_N})\right)^{\alpha_N} \\
&&\qquad\qquad\qquad \left(\frac{D_{x_1}}{2}-y_1+T_1(D_{y_1})\right)^{\beta_1}\dots \left(\frac{D_{x_N}}{2}-y_N+T_N(D_{y_N})\right)^{\beta_N}.
\end{eqnarray*}
Following the same procedure as in \cite{8} we obtain that the twisted Laplacian \eqref{TL}, as well as the operators in $\mathbb{R}^2$
$$
\left(x-\frac{1}{2}D_y+Q(D_x)\right)^2+\left(y+\frac{1}{2}D_x+R(D_y)\right)^2
$$
and
$$
(x-D_y+Q(D_x))^2+(y+R(D_y))^2,
$$
for arbitrary differential operators $Q(D_x)$ and $R(D_y)$ with real constant coefficients, are $\mathcal{S}_{\omega_1\oplus\omega_2}^{\sigma_1\oplus\sigma_2}$-regular, for every weight functions $\omega_1$, $\omega_2$, $\sigma_1$, $\sigma_2$.

\end{document}